\documentstyle[amstex]{amsart}

\makeatletter

\begin{document}
\bibliographystyle{amsalpha}

\newcommand{\e}{\epsilon}
\newcommand{\0}{{\bold 0}}
\newcommand{\w}{{\bold w}}
\newcommand{\y}{{\bold y}}
\newcommand{\balpha}{{\boldsymbol \alpha}}
\newcommand{\bt}{{\bold t}}
\newcommand{\z}{{\bold z}}
\newcommand{\x}{{\bold x}}
\newcommand{\N}{{\bold N}}
\newcommand{\Z}{{\bold Z}}
\newcommand{\F}{{\bold F}}
\newcommand{\R}{{\bold R}}
\newcommand{\Q}{{\bold Q}}
\newcommand{\C}{{\bold C}}
\newcommand{\BP}{{\bold P}}
\newcommand{\cO}{{\mathcal O}}
\newcommand{\cX}{{\mathcal X}}
\newcommand{\cH}{{\mathcal H}}
\newcommand{\cM}{{\mathcal M}}
\newcommand{\cD}{{\mathcal D}}
\newcommand{\cB}{{\mathcal B}}
\newcommand{\sB}{{\sf B}}
\newcommand{\cT}{{\mathcal T}}
\newcommand{\cI}{{\mathcal I}}
\newcommand{\cS}{{\mathcal S}}
\newcommand{\sE}{{\sf E}}

\newcommand{\sA}{{\sf A}}
\newcommand{\ga}{{\sf a}}
\newcommand{\es}{{\sf s}}
\newcommand{\m}{{\bold m}}
\newcommand{\bS}{{\bold S}}
\newcommand{\ovf}{{\overline{f}}}

\newcommand{\ihra}{\stackrel{i}{\hookrightarrow}}
\newcommand\rank{\mathop{\rm rank}\nolimits}
\newcommand\im{\mathop{\rm Im}\nolimits}
\newcommand\coker{\mathop{\rm coker}\nolimits}
\newcommand\Li{\mathop{\rm Li}\nolimits}
\newcommand\NS{\mathop{\rm NS}\nolimits}
\newcommand\Hom{\mathop{\rm Hom}\nolimits}
\newcommand\Ext{\mathop{\rm Ext}\nolimits}
\newcommand\Pic{\mathop{\rm Pic}\nolimits}
\newcommand\Spec{\mathop{\rm Spec}\nolimits}
\newcommand\Hilb{\mathop{\rm Hilb}\nolimits}
\newcommand\Ker{\mathop{\rm Ker}\nolimits}
\newcommand{\length}{\mathop{\rm length}\nolimits}
\newcommand{\res}{\mathop{\sf res}\nolimits}

\newcommand\lra{\longrightarrow}
\newcommand\ra{\rightarrow}
\newcommand\cJ{{\mathcal J}}
\newcommand\JG{J_{\Gamma}}
\newcommand{\wvskp}{\vspace{1cm}}
\newcommand{\vskp}{\vspace{5mm}}
\newcommand{\nvskp}{\vspace{1mm}}
\newcommand{\nid}{\noindent}
\newcommand{\new}{\nvskp \nid}
\newtheorem{Assumption}{Assumption}[section]
\newtheorem{Theorem}{Theorem}[section]
\newtheorem{Lemma}{Lemma}[section]
\newtheorem{Remark}{Remark}[section]
\newtheorem{Corollary}{Corollary}[section]
\newtheorem{Conjecture}{Conjecture}[section]
\newtheorem{Proposition}{Proposition}[section]
\newtheorem{Example}{Example}[section]
\newtheorem{Definition}{Definition}[section]
\newtheorem{Question}{Question}[section]
\renewcommand{\thesubsection}{\it}

\baselineskip=14pt
\addtocounter{section}{-1}

\title{Deformation of Okamoto--Painlev\'e Pairs and  
Painlev\'e equations}

\author{Masa-Hiko Saito,  Taro Takebe \& Hitomi Terajima}
\thanks{Partly supported by Grant-in Aid
for Scientific Research (B-09440015), (B-12440008) and (C-11874008), 
the Ministry of
Education, Science and Culture, Japan }
\address{Department of Mathematics, Faculty of Science, 
Kobe University, Kobe, Rokko, 657-8501, Japan}
\email{mhsaito@math.kobe-u.ac.jp}
\email{takebe@math.kobe-u.ac.jp}
\email{terajima@math.kobe-u.ac.jp}
\date{July, 15, 2000}
\begin{abstract}
In this paper, we introduce the notion of generalized 
rational Okamoto--Painlev\'e pair $(S, Y)$ by generalizing  the notion of the spaces of initial conditions of Painlev\'e equations.  After classifying those pairs, 
we will establish an algebro-geometric approach to derive the Painlev\'e
differential equations from the deformation of Okamoto--Painlev\'e pairs by using the  local cohomology groups.  
Moreover the reason why the 
Painlev\'e equations can be written in Hamiltonian systems 
is clarified by means of the holomorphic symplectic 
structure on $S - Y$. Hamiltonian structures for Okamoto--Painlev\'e pairs of type   $\tilde{E}_7 (= P_{II})$ and $\tilde{D}_8 (= P_{III}^{\tilde{D}_8})$ are calculated explicitly as examples of our theory.  
\end{abstract}
\maketitle

\section{Introduction}

In the study of Painlev\'e equations, the spaces of initial conditions introduced by K. Okamoto \cite{O1}, \cite{O2}, \cite{O3} have been playing essential roles.    
It is known that each 
Painlev\'e differential equation  
 is equivalent to one of  
Hamiltonian systems whose Hamiltonians are given by the polynomials in two variables $(x, y)$. (See Table \ref{tab:painleve} and  \ref{tab:ham} in \S \ref{sec:painleve}). 
The space 
$(x, y) \in \C^2 $ can be compactified and one can obtain 
a pair $ (S, Y) $ of complex projective surface $S$ and an anti-canonical divisor $ Y \in |-K_S| $ 
such that $S - Y_{red}$ becomes a space of initial conditions.  
  In the study of the space of initial conditions as in \cite{O1}, \cite{MMT}, it became  clear that  after 
  eliminating the singularities of Painlev\'e differential 
  equation by blowings-up,   the boundary divisor  $Y$ should have  the same configuration as in the list of degenerate 
  elliptic curves classified by Kodaira \cite{Kod}.  
 This condition can be translated into the following conditions.  Let $Y= \sum_{i=1}^r m_i Y_i \in |-K_S|  $ 
 be the irreducible decomposition.  Then $Y$ is called 
{\em of canonical type} if and only if
 $$
 \deg (- K_S)_{|Y_i} = \deg Y_{|Y_i} = Y \cdot Y_i = 0 \quad \mbox{for all $ i $. }
 $$  

 In \cite{Sa-Tak}, we call such a pair $ (S, Y ) $ an {\em Okamoto--Painlev\'e pair} if  $ S - Y_{red} $ contains $\C^2$ as a Zariski open set and $F = S -\C^2$ is a normal crossing divisor.  One can verify that 
all compactifications of the spaces of initial conditions 
of known Painlev\'e equations satisfy these conditions (cf. \cite{O1}, \cite{MMT}). Therefore,  in this notation, the former studies of Painlev\'e equations establish the route in the direction:
 \begin{equation}
 \fbox{Painlev\'e equations} \quad  \Longrightarrow \quad  \fbox{Okamoto--Painlev\'e pairs $(S, Y)$} . 
 \end{equation}
 
The main purpose of this paper is to establish the route 
backward, that is, the route in the following direction  
\begin{equation}
\fbox{Okamoto--Painlev\'e pairs $(S, Y)$} 
\quad \Longrightarrow \quad  \fbox{Painlev\'e equations}. 
\end{equation}

 Our main tool here is the deformation theory of pairs $(S, Y)$ (cf. \cite{KS}, \cite{Kaw}) and local cohomology exact sequence (cf. \cite{B-W}).    
 The deformation theory was  established  by Kodaira--Spencer in \cite{KS} in late 1950's.  Kawamata \cite{Kaw} generalized the deformation theory to compactified complex manifolds, or  pairs of smooth compact complex manifolds and simple normal crossing  subvarieties.  Applying the deformation theory, we can see that the space of infinitesimal deformations of the Okamoto--Painlev\'e pair $(S, Y)$ is isomorphic to the cohomology group $H^1(S, \Theta_{S}(- \log D))$ where $D = Y_{red}$. Looking at the restriction homomorphism
$$
\mbox{\sf res}: H^1(S, \Theta_{S}(- \log D)) \lra H^1(S-D, \Theta_{S-D}), 
$$   
 we see that the kernel of the restriction map $\mbox{\sf res}$ has the important meaning,  that is,  
 the direction corresponding to the kernel of $\mbox{\sf res}$ 
 is the infinitesimal deformation of $(S, Y)$ which induces the trivial deformation on $ S - D$. 
Roughly speaking, one can say that the Painlev\'e differential equations describe the deformations corresponding to the direction of the kernel of the restriction map.  

To be precise, let us consider a family of Okamoto--Painlev\'e pairs  
$\cD \hookrightarrow \cS  \lra \cB_R$ with one-dimensional  base space $\cB_R$ with a coordinate $t$ such that the Kodaira--Spencer class $\rho(\frac{\partial}{\partial t}) $ lies in the kernel of the restriction map $\mbox{\sf res}$.  Then by using the affine covering of $\cS - \cD$  and \v{C}ech cocycles, we can derive a system of ordinary differential equation.  Note that this observation will be applicable to the higher dimensional cases.   

From these observation, we see that  the kernel of $\mbox{\sf res}$ corresponds to the directions of  time variables in the Painlev\'e differential equation.  
  
Furthermore, we can apply the local cohomology exact sequence to our settings and obtain the exact sequence (cf. \cite{B-W}, \cite{Gr})
$$\
H^1_{D}(S, \Theta_{S}( - \log D)) \stackrel{\mu}{\lra} H^1(S, \Theta_{S}(- \log D)) \stackrel{\mbox{\sf res}}{\lra} H^1(S-D, \Theta_{S-D}).  
$$  
Under the  condition that $(S,Y)$ is of non-fibered type, 
the map $\mu$ is injective, and hence, the local cohomology group $H^1_{D}(S, \Theta_{S}( - \log D))$ coincides 
with the kernel of $\mbox{\sf res}$.  
Therefore, non-zero element of 
the  local cohomology group $H^1_{D}(S, \Theta_{S}( - \log D))$ corresponds to a time variable of the Painlev\'e equation.  
For a generalized rational Okamoto--Painlev\'e pair $(S, Y)$ of additive type, if $Y_{red}$ is normal crossing divisor,  Terajima \cite{T} proved that 
the dimension of the local cohomology group  is positive,  
hence,  we can always obtain a differential equation.   

In order to obtain an explicit 
differential equation for each type $R$ from our 
setting, we need to construct a global family of 
generalized rational Okamoto--Painlev\'e pairs of type 
$R$  over a parameter space $\cM_R \times \cB_R$
which is semiuniversal at each point.  Moreover 
we need to  introduce a good affine open covering of 
the total space such that the rational two form $\omega_S$ 
restricted to $S - Y_{red}$ has a canonical form.  
 (See \S \ref{sec:globaldef} and \S \ref{sec:globaltoham}.) 

The organization of this  paper is  as follows.  In \S \ref{genop}, we define the notion of generalized Okamoto--Painlev\'e pairs and recall 
the relation to generalized Halphen surfaces, which are 
studied by Sakai \cite{Sakai}.  We also classify  generalized rational Okamoto--Painlev\'e pairs $(S, Y)$ such that $Y_{red}$ are normal crossing divisors.  A generalized rational Okamoto--Painlev\'e pair $(S, Y)$ is called of fibered type if there exists a structure of elliptic fibration $f:S \lra \BP^1$ such that $f^*(\infty) = n Y$ for some $n \geq 1$.   We show that $(S, Y)$ is not of fibered type if and only
 if regular algebraic functions on  $S - Y_{red}$ are just constant functions. This fact is also important for later purpose. After recalling the theory of deformation of pairs in \S \ref{sec:deformation}, we investigate the cohomology groups for generalized rational Okamoto--Painlev\'e pairs.  In \S 3, we will apply 
the theory of local cohomology to our situation, and 
have the fundamental exact sequence (Proposition \ref{prop:local}).  Moreover, we state an important result,  Theorem \ref{thm:time}, which is proved in \cite{T}.   
After reviewing the Kodaira--Spencer theory in \S \ref{sec:KS}, in \S \ref{sec:globaldef},   we will explain how one can construct the family of generalized 
rational Okamoto--Painlev\'e pairs and their open coverings. In \S \ref{sec:globaltoham}, we will state our main theorem (Theorem 6.1), which states that from the special 
global deformation of generalized 
Okamoto--Painlev\'e pairs one can obtain the differential 
equations.  Moreover the reason for the equations to be 
in Hamiltonian systems will be explained geometrically.  
In \S \ref{sec:examples}, 
we will derive the Painlev\'e equations from 
the deformations of Okamoto--Painlev\'e pairs of type 
$\tilde{E}_7$ and $\tilde{D}_8$.

 Prior to our work here, in \cite{SU},  
M.-H. Saito and H.  Umemura 
essentially pointed out that the deformation of spaces of 
initial conditions describes Painlev\'e equation completely.   In this sense, this paper is a continuation of \cite{SU}, though  we have clarified the meaning of time variables by means of local cohomology groups in this paper.  

The recent work due to Sakai \cite{Sakai} 
introduce  the following beautiful viewpoint: The geometry of certain rational surfaces with the symmetries  
induced by Cremona transformations describe the  discrete Painlev\'e equations and the Painlev\'e equation can be obtained as a limit of the discrete Painlev\'e equations.  We owe much to his beautiful paper \cite{Sakai}.  In particular, some of the explicit calculations are done by using his descriptions of the 
family of Okamoto--Painlev\'e pairs in \cite{Sakai}.  

The works of Takano et al \cite{MMT}, \cite{ST} is also essential to our work.  In \S \ref{sec:examples}, we use the coordinate systems introduced by them.  

  The series of the work is started by \cite{Sa-Tak}, where 
  we introduce the notion of Okamoto-Painlev\'e pair and 
  classify  Okamoto--Painlev\'e pairs $(S, Y)$.  
  
One of  missing points in our work here is the theory of B\"{a}cklund transformation of Okamoto--Painlev\'e pairs.  In this direction, one should refer to  a series of works of M. Noumi and Y. Yamada (e.g. \cite{NY}), also Sakai's work \cite{Sakai}.  In \cite{SU}, the authors tried to understand the B\"acklund transformation by using the notion of ``flip " or ``flop" in 
the theory of the minimal models of higher dimensional varieties. We will investigate this point in future.

\section{Generalized Okamoto--Painlev\'e Pairs}
\label{genop}

Let $S$ be  a complex projective 
surface. We denote by $K_S$ the canonical 
line bundle or the canonical divisor 
class  of $S$.  Assume that the anti-canonical divisor class  $-K_S$ is effective, that is, there exists an effective 
divisor $Y \in |-K_S|$. Geometrically, this is equivalent to the existence of  a  rational 2-form  $\omega_Y$ on $S$ whose 
corresponding  divisor $(\omega_Y)=-(\omega_Y)_{\infty} = -Y$. Such a divisor $Y$ is called an anti-canonical divisor of $S$ as usual. Since $\omega_Y$ dose not vanish on $S - Y$, it induces a 
holomorphic symplectic structure  on $S-Y$.

In \cite{Sa-Tak}, we introduce a notion of Okamoto--Painlev\'e pair $(S, Y)$ which is a pair of  complex projective  surface $S$ and an 
anti-canonical divisor $Y$ satisfying certain conditions ([Definition 2.1 \cite{Sa-Tak}]).  Generalizing the notion, we will start this section with the following definition.  

\begin{Definition}\label{def:op}
{\rm
Let $(S, Y)$ be a pair  of a  complex projective surface $S$  and an anti-canonical divisor $Y \in |-K_S|$ of $S$. Let $Y= \sum_{i=1}^r m_i Y_i$ be the irreducible 
decomposition of $Y$.    
 We call a pair $(S, Y)$ a {\em generalized  Okamoto--Painlev\'e Pair} if for all $i, 1 \leq i \leq r$, 
\begin{equation} \label{eq:canonical}
Y \cdot Y_i = \deg Y_{|Y_i} = 0 .
\end{equation}
}
\end{Definition}

According to  \cite{Sa-Tak},  we listed 
the additional conditions for Okamoto--Painlev\'e pairs besides the condition  (\ref{eq:canonical})   as follows.  
\begin{enumerate}
\item Let us set  $ D :=  Y_{red} = \sum_{i=1}^{r} Y_i $.
Then $ S - D $ contains the complex affine plane $ \C^2$ as a Zariski open set.

\item Set $ F = S - \C^2 $ where $ \C^2 $ is the same Zariski open set  as in (1). 
Then $ F $ is a (reduced) divisor with normal crossings. In particular, $ D=Y_{red} $ is also a reduced divisor with normal crossings. 
\end{enumerate}

Under this definition, we proved the following 
classification of Okamoto--Painlev\'e pairs in \cite{Sa-Tak}. We remark that an Okamoto--Painlev\'e pair of type $\tilde{D}_7$ did not appear in the classification of classical Painlev\'e equations  \cite{O1}.

\begin{Theorem}\label{thm:classf}
 $($\cite{Sa-Tak}.$)$
Let $(S, Y)$ be a generalized Okamoto-Painlev\'e pair and assume that $S - Y_{red}$ contains $\C^2$ and $F = S-\C^2$ is a reduced divisor with normal crossings. $($That is, $(S,Y)$ is an Okamoto--Painlev\'e pair in original sense. $)$ 
Then we have the following assertions.

\begin{enumerate}

\item The surface $S$ is a projective rational surface.

\item  The configuration of $Y$ counting with
 multiplicity  is in the list of 
 Kodaira's classification of singular fibers of 
 elliptic surfaces $($cf. \cite{Kod}$)$.   More precisely, 
it coincides with  one in  the following Table \ref{tab:op}.  $($In Figure \ref{fig:dynkin}, each 
 line denotes a smooth rational curve $C$ with $C^2= -2$ 
 and the configuration of lines show  how they  
 intersect to each other. The number next to each line 
 denotes the multiplicity of each curve in $Y = -K_S$.$)$    
 
 \end{enumerate}
\end{Theorem}

\vspace{0.5cm}

\begin{table}[h]
\begin{center}
\begin{tabular}{||c||c|c|c|c|c|c|c|c||} \hline
    &   & & & & & &  \\
$Y$ & $\tilde{E_8}$ & $\tilde{E_7}$ & $\tilde{D_7}$ & 
$\tilde{D_6}$ & $\tilde{E_6}$ &$\tilde{D_5}$ & $\tilde{D_4}$ \\
    &   & & & & & &  \\ \hline
        &   & & & & & &  \\
Kodaira's notation & $II^*$ & $III^*$ & $I_{3}^*$ & $I_{2}^*$ & 
$IV^*$ & $I_{1}^*$ & $I_{0}^*$ \\ 
    &   & & & & & &  \\ \hline
        &   & & & & & &  \\ 
Painlev\'{e} equation & $P_{I}$  & $P_{II}$ & $P^{\tilde{D}_7}_{III}$ & $P_{III} = P^{\tilde{D}_6}_{III}$ & $P_{IV}$
& $P_{V}$  &  $P_{VI}$ \\
    &   & & & & & &  \\ \hline
\end{tabular}
\vspace{0.5cm}
\caption{}
\label{tab:op}
\end{center}
\end{table}

\begin{figure}
\vspace{0.5cm}
\begin{center}
\begin{picture}(150,150)(140,0)
\put(0, -2){\line(0,1){48}}
\put(-4, 40){\line(1,0){48}}
\put(40, 36 ){\line(0,1){48}}
\put(36, 80){\line(1,0){48}}
\put(80, 76){\line(0,1){48}}
\put(76, 120){\line(1,0){88}}
\put(120, 112){\line(0,1){48}}
\put(160,124){\line(0, -1){48}}
\put(156, 80){\line(1,0){48}}
\put(-10, 18){$1$}
\put(20, 46){$2$}
\put(32, 58){$3$}
\put(60, 86){$4$}
\put(72, 98){$5$}
\put(100, 125){$6$}
\put(125, 140){$3$}
\put(165, 98){$4$}
\put(185, 85 ){$2$}
\put(100, 10){$\tilde{E_8}$}
\put(250, 36){\line(0,1){48}}
\put(242, 60){$1$}
\put(246, 80){\line(1,0){48}}
\put(270,85){$2$}
\put(290, 76){\line(0,1){48}}
\put(280,100){$3$}
\put(286, 120){\line(1,0){88}}
\put(310, 125){$4$}
\put(330, 116){\line(0,1){48}}
\put(322, 140){$2$} 
\put(370, 124){\line(0, -1){48}}
\put(375, 100){$3$}
\put(368, 80){\line(1,0){48}}
\put(390, 85){$2$}
\put(410, 84){\line(0, -1){48}}
\put(415, 60){$1$}
\put(328,10){$\tilde{E_7}$}
\end{picture}

\vspace{0.5cm}

\begin{picture}(150,150)(140,30)
\put(-4, 80){\line(1,0){48}}
\put(40, 76){\line(0,1){48}}
\put(36, 120){\line(1,0){88}}
\put(80, 116){\line(0,1){48}}
\put(56,150){\line(1,0){48}}
\put(60, 160){$1$}
\put(120,124){\line(0, -1){48}}
\put(116, 80){\line(1,0){48}}
\put(20, 86){$1$}
\put(32, 98){$2$}
\put(60, 125){$3$}
\put(85, 130){$2$}
\put(125, 98){$2$}
\put(145, 85 ){$1$}
\put(80, 50){$\tilde{E_6}$}
\put(200,80){\line(1,1){60}}
\put(220,90){$1$}
\put(210,70){\line(1,1){60}}
\put(230,80){$1$}
\put(220, 140){\line(1,-1){64}}
\put(252, 90){$2$}
\put(277, 77){\line(1,1){64}}
\put(300, 90){$2$}
\put(300,140){\line(1,-1){64}}
\put(330, 90){$2$}
\put(357, 77){\line(1,1){64}}
\put(380, 90){$2$}
\put(440,80){\line(-1,1){64}}
\put(415, 90){$1$}
\put(430, 70){\line(-1,1){64}}
\put(405,80){$1$}
\put(300, 50){$\tilde{D_7}$}
\end{picture}

\begin{picture}(150,150)(140,0)
\put(0,80){\line(1,1){60}}
\put(20,90){$1$}
\put(10,70){\line(1,1){60}}
\put(30,80){$1$}
\put(20, 140){\line(1,-1){64}}
\put(52, 90){$2$}
\put(108, 90){$2$}
\put(68, 85){\line(1,0){88}}
\put(137, 77){\line(1,1){64}}
\put(160, 90){$2$}
\put(210,80){\line(-1,1){64}}
\put(185, 90){$1$}
\put(200, 70){\line(-1,1){64}}
\put(175,80){$1$}
\put(100, 20){$\tilde{D_6}$}
\put(260,80){\line(1,1){50}}
\put(280,90){$1$}
\put(270,70){\line(1,1){50}}
\put(290,80){$1$}
\put(280, 140){\line(1,-1){64}}
\put(312, 90){$2$}
\put(317, 77){\line(1,1){64}}
\put(337, 90){$2$}
\put(400,80){\line(-1,1){50}}
\put(380,85){$1$}
\put(390,70){\line(-1,1){50}}
\put(370,75){$1$}
\put(320, 20){$\tilde{D_5}$}
\end{picture}

\begin{picture}(100, 100)(160,0)
\put(10,80){\line(1,0){140}}
\put(80,85){$2$}
\put(25,90){\line(0,-1){50}}
\put(18,65){$1$}
\put(40, 90){\line(0,-1){50}}
\put(33,65){$1$}
\put(120,90){\line(0,-1){50}}
\put(123,65){$1$}
\put(135,90){\line(0,-1){50}}
\put(138,65){$1$}
\put(80, 20){$\tilde{D_4}$}
\end{picture}
\end{center}
\caption{}
\label{fig:dynkin}
\end{figure}

\noindent
{\bf Generalized Halphen surfaces}

\par
\vspace{0.5cm}
According to Sakai [\S 4, \cite{Sakai}], we recall the following definition.  

\begin{Definition}{\rm 
\begin{enumerate}

\item  Let $S$ be a rational surface with an effective anti-canonical  divisor $Y \in | -K_S|$.  Let $Y = \sum_{i=1}^r m_i Y_i$ be the irreducible decomposition.  The divisor  $Y$ is called {\em of canonical type} if
$$
K_S \cdot Y_i = - Y \cdot Y_i = 0 \quad \mbox{for all}  \quad i 
$$

\item A rational surface $S$ is called a generalized Halphen surface if $S$ has  an effective anti-canonical divisor $Y$ of canonical type. A generalized Halphen surface $S$ is called of index one if
$$
 \dim |-K_S| = 1.
$$  
\end{enumerate} }
\end{Definition}

\begin{Remark} {\rm By Riemann-Roch theorem, it is easy to see that for a generalized Halphen surface $S$ $\dim |-K_S| \leq 1$. } 
\end{Remark}

The following Proposition ensures that  one can obtain 
a generalized Halphen surfaces from  blowing-up of 9-points of  
$\BP^2$.  

\begin{Proposition}\label{prop:blowup}
 $($ Proposition 2, \S 2, \cite{Sakai} $)$.    Let $S$ be a generalized Halphen surface, then there exists a birational morphism $ \rho:S \lra \BP^2$.   
\end{Proposition}

Let $(S, Y)$ be a generalized Okamoto-Painlev\'e pair such that $S$ is a rational surface.  
Then $S$ is a generalized 
Halphen surface with a specified anti-canonical divisor $Y$. As a corollary of Proposition \ref{prop:blowup}, we obtain the following corollary.

\begin{Corollary} Let $(S, Y)$ be a generalized rational Okamoto--Painlev\'e pair.  Then $S$ can be obtained as 9 points blowing-up of $\BP^2$. 
\end{Corollary}

 One can show that $Y$ has a same configuration 
as one of Kodaira's degenerate elliptic curves for a 
generalized rational Okamoto--Painlev\'e pair $(S, Y)$ (cf. Proof of Theorem 2.1 in \cite{Sa-Tak}).

    If $S$ is 
a generalized Halphen surface of index one,  the morphism associated to 
the linear system $|-K_S|$ induces an elliptic fibration 
$\varphi:S \lra \BP^1$ with $\varphi^*(\infty) = Y$. 
(Here $\varphi^*(\infty)$ denotes the scheme theoretic 
fibers at $\infty \in \BP^1$. ) 
This leads us the following 
\begin{Definition}{\rm 
A generalized Okamoto-Painlev\'e pair $(S, Y)$ is called 
``of fibered  type" if there exists an elliptic fibration 
$\varphi:S \lra \BP^1$ such that 
$\varphi^*(\infty) = n Y$ for some $n \geq 1$.   If $(S,Y)$ is not of fibered type, we call $(S, Y)$ ``{\em of non-fibered type} ".  
}
\end{Definition}

Note that if $(S, Y)$ is of fibered type and 
$ \varphi:S \lra \BP^1 $ is elliptic surface with 
$\varphi^*(\infty) = n Y$ with $n>1$,  $\varphi^*(\infty)$ is called a multiple fiber.  This happens only when  $Y$ is  
of elliptic type or multiplicative type in the notation below.  

In \cite{Sakai}, Sakai classified generalized Halphen surface $S$ with $\dim |-K_S| = 0$.   
In the case when $\dim |-K_S| = 0$,  the associated Okamoto-Painlev\'e pair $(S, Y)$ with a unique member $Y \in |-K_S|$
is of non-fibered type and they 
 can be classified by means of the configuration of $Y$.

Let $Y = \sum_{i=1}^r m_i Y_i$ be the irreducible 
decomposition of $Y$.  
Denote by $M(Y)$ the sublattice of $\Pic(S) \simeq H^2(S, \Z)$ 
generated by the irreducible components $\{ Y_i \}_{i=1}^r$.  Here the bilinear form on $\Pic(S)$ is $(-1)$ times 
the intersection form on $\Pic(S)$.  Then $\{ Y_i \}_{i=1}^r $ forms a root basis of $M(Y)$ and we denote by $R(Y)$ the type of the root system. 

One can easily classify  $R(Y)$ as in Table 
\ref{tab:genop}.  
Note that according to the type of $Y$,  $R(Y)$ can be classified into three classes: {\it elliptic type} when $Y$ is a smooth 
elliptic curve, {\it multiplicative type} when $Y$ is a cycle of 
rational curves, {\it additive type}  when the configuration of $Y$ is tree.  These types also correspond to the type of generalized Jacobean $\Pic^0(Y)$ of $Y$.  If we denote by 
$(\Pic^0(Y))^0 $ the component of identity of $\Pic^0(Y)$, 
we have the following correspondence (cf. Table \ref{tab:pic}).  

\vspace{0.5cm}

\begin{table}[h]
\begin{center}
\begin{tabular}{|l|l|} \hline
 $R(Y)$        &    $(\Pic^0(Y))^0 $  \\ \hline  
elliptic type &  smooth elliptic curve $Y$  \\ \hline 
multiplicative type  & ${\bf G}_m \simeq \C^{\times}$ \\ \hline 
additive  type &  ${\bf G}_a \simeq \C$  \\ \hline
\end{tabular}
\vspace{0.2cm}
\end{center}
\caption{}
\label{tab:pic}
\end{table}

\begin{table}[h]
\begin{center}
\begin{tabular}{|l|l|}
\hline 
   &     \\
  & $R(Y)$ (Kodaira type) \\ 
  &    \\
  \hline 
    &  \\ 
\mbox{elliptic type}  & $ \tilde{A_0} (= I_0)$  \\
   &  \\ \hline  
                      &     \\
\mbox{multiplicative type} & $\tilde{A_0}^*(=I_1), 
\tilde{A_1}(=I_2), \cdots, \tilde{A_7} (= I_8), \tilde{A_8} (=I_9)$  \\ 
   &  \\ 
\hline 
     &     \\
\mbox{additive type} & $\tilde{A_{0}^{**}} (=II), \tilde{A_1^{*}} (= III), \tilde{A_2^{*}}( =IV) $  \\ 
    &    \\
     & $\tilde{D_4}(=I_0^*), \cdots,  \tilde{D_8}(=I_4)$ \\ 
     &   \\
     & $ \tilde{E_6}(=IV^*),          \tilde{E_7}(=III^*), \tilde{E_8} (= II^*)$ \\
        &  \\ 
     \hline 
\end{tabular}
\vspace{0.3cm}
\caption{}
\label{tab:genop}
\end{center}
\end{table}

\begin{Proposition}\label{prop:classf}
Let $(S, Y)$ be a generalized rational Okamoto--Painlev\'e 
pair such that $Y_{red}$ is a divisor with only normal crossings.    Then besides the list of Okamoto--Painlev\'e pairs in Table 1, we have a pair $(S, Y)$ of type $\tilde{D_8}$ and also $\tilde{A_{r}}$ for $0 \leq r \leq 8$ and $\tilde{A_0}^*$.  Here for $\tilde{A_0}$, $Y$ is a smooth elliptic curve (Kodaira $I_0$-type)  and for $\tilde{A_0}^*$,  $Y$ is a rational curve with a node (Kodaira $I_1$-type). 

\end{Proposition}

We list up  generalized rational Okamoto--Painlev\'e  pairs 
with normal crossing divisor $ Y_{red}$ in Table \ref{tab:addop}.

\vspace{0.5cm}
\begin{table}[h]
\begin{center}
\begin{tabular}{||c||c|c|c|c|c|c|c|c|c|c|c||} \hline
    &  &  & & & & & & &  & &  \\
$Y$ & $\tilde{E_8}$ & $\tilde{D_8}$ & $\tilde{E_7}$ & $\tilde{D_7}$ & 
$\tilde{D_6}$ & $\tilde{E_6}$ &$ \tilde{D_5} $ & $ \tilde{D_4} $ & $ \tilde{A}_{r-1} $ & $\tilde{A_0} $  & $\tilde{A_0}^* $   \\
    & &   & & & & & &  & {\tiny $2 \leq r \leq 9$} & $ r=1$ & $r =1 $   \\ \hline
  && &  &  &   &   &   &  & & & \\
Kodaira's notation & $II^*$ & $I_{4}^*$  & $III^*$ & $I_{3}^*$ & $I_{2}^*$ & 
$IV^*$ & $I_{1}^*$ & $I_{0}^*$ & $I_{r}$ & $I_0$ & $I_1 $  \\ 
    &  & & & & & & &  & & &\\ \hline
        & &   & & & & & &  & & & \\ 
$r= \sharp$ of comps. of $Y$  & 9   & 9  & 8 &8 &7 &7 &6 &5  & $r$ & 1 & 1 \\ 
     &   & & & & & & &  & & &    \\ \hline 
\end{tabular}
\vspace{0.3cm}
\caption{}
\label{tab:addop}
\end{center}
\end{table}

\begin{figure}[h]
\begin{center}
\begin{picture}(150,150)(270, 30)
\put(200,80){\line(1,1){60}}
\put(220,90){$1$}
\put(210,70){\line(1,1){60}}
\put(230,80){$1$}
\put(220, 140){\line(1,-1){64}}
\put(252, 90){$2$}
\put(277, 77){\line(1,1){64}}
\put(300, 90){$2$}
\put(300, 120){\line(1,0){90}}
\put(340, 110){$2$}
\put(350,140){\line(1,-1){64}}
\put(380, 90){$2$}
\put(407, 77){\line(1,1){64}}
\put(430, 90){$2$}
\put(490,80){\line(-1,1){64}}
\put(465, 90){$1$}
\put(480, 70){\line(-1,1){64}}
\put(455,80){$1$}
\put(330, 50){$\tilde{D_8}$}
\end{picture}

\begin{picture}(400,250)(20, 0)
\put(150,50){\line(1,0){130}}
\put(210,55){$1$}
\put(250,40){\line(1,1){65}}
\put(270,70){$1$}
\put(300,80){\line(1,4){20}}
\put(300,120){$1$}
\put(320,140){\line(-1,2){40}}
\put(295,170){$1$}
\put(300,200){\line(-2,1){100}}
\put(250,210){$1$}
\put(180,40){\line(-1,1){65}}
\put(153,70){$1$}
\put(130,80){\line(-1,4){20}}
\put(123,120){$1$}
\put(110,140){\line(1,2){40}}
\put(132,170){$1$}
\put(130,200){\line(2,1){100}}
\put(210, 20){$\tilde{A_8}$}
\put(170,210){$1$}
\end{picture}
\end{center}
\vspace{0.3cm}
\caption{}
\label{fig:d8a8}

\end{figure}

\noindent
{\bf Regular algebraic functions on $S - Y$.}

\vspace{0.5cm}
Let $(S, Y)$ be a generalized rational Okamoto-Painlev\'e pair.  
If $(S, Y)$ is of fibered type, that is, if  there exists an elliptic fibration $\varphi:S \lra \BP^1$ with $\varphi(\infty)^* = nY$, pulling back  non-constant regular algebraic functions on $\BP^1 - \{\infty \} \simeq \C$, 
we have many  non-constant regular functions on $S - Y$.  We can prove the 
converse of this fact.
\begin{Proposition}\label{prop:reg}
  Let $(S, Y)$ be a generalized rational Okamoto--Painlev\'e pair.  
  The following conditions are equivalent to each other.
\begin{enumerate}  
 \item $(S, Y)$ is of non-fibered type.
 \item $H^0(S -Y, \cO^{alg}) \simeq \C$, that is, all regular algebraic functions of $S-Y$ are constant functions. 
\end{enumerate}
\end{Proposition}
 
{\it Proof.}
As we remarked as  above,  the implication $(2) \Rightarrow (1) $ is obvious.  
Assume that there exists a non-constant regular function $f$ on $S-Y$.  
Then  the morphism $f:S-Y \lra \C$ extends to a morphism 
$$
\overline{f}:S \lra 
\BP^1.
$$  
Set $Y' = \ovf^*(\infty)$.  Since $\ovf$ is regular on $S - Y$, recalling 
the irreducible decomposition of $Y = \sum_{i=1}^r m_i Y_i$, we can 
write 
$$
Y' = \sum_{i=1}^r a_i Y_i
$$
with $a_i \geq 0$.  First we show that $a_i > 0$ for every $1 \leq i \leq r$.  
If  $a_i = 0$ for some $i$, the configuration of $Y'_{red}$ becomes a proper 
sub-graph of the configuration of $Y_{red}$.  Since the graph of $Y_{red}$ corresponds to an affine Dynkin diagram, one can easily see that 
$Y'$ can be contracted to rational double points $\{p_1, \cdots, p_s \}$ and obtain a normal surface $S'$ with normal singular points  $\{p_1, \cdots, p_s \}$. 
Since $S - Y'$ and $S' - \{p_1, \cdots, p_s \}$ are isomorphic and $\ovf$ is regular on $S-Y' \simeq S' - \{p_1, \cdots, p_s \}$, $\ovf$ extends to a regular function on $S'$. Since $S'$ is proper, $\ovf$ must be constant and hence $f$ is also constant.  This contradicts to the original assumption.  
Hence we see that  $a_i > 0$ for all $1 \leq i \leq r$.   

This implies that $Y'_{red} = Y_{red}$ and since $Y_{red}$ is connected, 
so is $Y'_{red}$.  
Taking the Stein factorization if necessary, we may assume that all of the 
fiber of $\ovf:S \lra \BP^1$ is connected and $\ovf^*(\infty)_{red} = Y'_{red} = Y_{red} $.  (Note that $S$ is a rational surface, hence the irregularity of $S$ is 
zero.)  We will show that general fiber of $\ovf$ is an elliptic curve. Since 
$S$ is smooth and $\ovf$ has connected fibers, we only have to show that the 
virtual genus of $Y'$ is one.  Since $-K_S = Y$ and $Y' = \sum_{i=1}^r a_i Y_i$, we see that 
$$
K_S \cdot Y' = \sum_{i=1}^r a_i (-Y) \cdot Y_i = 0
$$
by  definition of Okamoto--Painlev\'e pair.  Moreover since $Y'$ is linear equivalent to a general fiber of $\ovf$, we see that  $Y' \cdot Y_i = 0$ for all $1 \leq i \leq r$.  Hence we see that $(Y')^2 = 0$.  Then the virtual genus of 
$Y'$ is given by 
$$
\pi(Y') =  \frac{K_S \cdot Y' + (Y')^2}{2}+1  = 1, 
$$
and this completes the proof of proposition.  
\qed

\section{Deformation of generalized rational Okamoto--Painlev\'e pairs}\label{sec:deformation}

\vspace{0.5cm}

In this section, we will recall necessary background of theory of infinitesimal deformation of Okamoto--Painlev\'{e} pairs. 
First, we will recall the general theory of deformation of pairs.

\vspace{0.5cm}
\noindent
{\bf  General Theory of Deformation of Pairs}

\vspace{0.5cm}

Let $(X, H) $ be a pair of a complex manifold $X$ and a 
(reduced) normal crossing divisor and let $H = \sum_{i=1}^r H_i$ be the irreducible decomposition of $H$. By a technical reason we will assume that $H$ is a simple normal crossing 
divisor, that is, each irreducible component $H_i$  is a smooth divisor.  We call such a pair $(S, H)$ is a {\em non-singular pair }.  

For such a non-singular pair $(S, H)$, 
the normalization $\tilde{H}$ of $H$ is given by the disjoint union $\coprod_{i=1}^r H_i$ of $H_i$'s, and we denote by $\nu:\tilde{H} \lra H$ the normalization map.  

First, we recall the general theory of deformation of 
a non-singular pair $(X, H)$ due to Kawamata \cite{Kaw}. (See also \cite{SSU}). 

Let $\Omega_X^1( \log H)$ denote the sheaf of germs of 
meromorphic one forms on $X$ 
which have logarithmic poles along $H$. Moreover, we set
$$
\Theta_{X}( - \log H ) := \underline{\Hom}( \Omega_{X}^1(\log H), \cO_X).
$$
This is the sheaf of germs of regular vector fields 
which have logarithmic zero along $H$.

\begin{Definition}  {\rm (Cf. [Definition 3, \cite{Kaw}])   }
{\rm A deformation of a non-singular pair $(X, H)$ is a 
5-tuple $({\cal X}, {\cal H}, \pi, B, \iota ) $ 
\begin{enumerate}
\item $ \pi:{\cal X} \lra B$ is a proper smooth 
holomorphic map  from a complex manifold ${\cal X}$ to a connected complex manifold $B$  

\item ${\cal H} = \sum_{i=1}^r {\cal H}_i$ is a simple 
normal crossing divisor of ${\cal X}$

\item For a point $0 \in B$, we have an isomorphism 
$ \iota: (\pi^{-1}(0), \pi^{-1}(0) \cap {\cal H}) = ({\cal X}_0, {\cal H}_0) \simeq (X, H)$.

\item $\pi$ is locally a projection of a product space as well as the restriction of it to ${\cal H}$, that is, for 
each $p \in {\cal X}$ there exists an open neighborhood $U$ of $p$ and an isomorphism $\varphi: U \lra V \times W$, where $V = \pi(U)$ and $W = \pi^{-1}(\pi(p))$, such that the following diagram 
$$
\begin{array}{ccc}
U & \stackrel{\varphi}{\lra} & V \times W \\
\pi \searrow &   & \swarrow  pr_1 \\
   &  V &   
\end{array} 
$$ 
 is commutative and $\varphi(U \cap {\cal H}) = V \times (W \cap {\cal H})$. 

The deformation of a pair 
is often denoted by the following diagram:
\end{enumerate}
$$
\begin{array}{ccl}
  {\cal X} &  \hookleftarrow & {\cal H}      \\
   \pi \downarrow \hspace{0.3cm} &    \swarrow &    \\
  B  &  & 
\end{array} 
$$

}
\end{Definition}

For a deformation $\pi:{\cal X} \lra B$ of complex manifold $X = {\cal X}_0$, we can define the Kodaira--Spencer class 
$$
\rho:T_0(B) \lra H^1(X, \Theta_{X} ).
$$

Similarly, for a deformation   of a pair $(X, H)$ as above, 
 we can define   the Kodaira--Spencer map
$$
\rho: T_0(B) \lra H^1(X, \Theta_X(- \log H)). 
$$
As for the existence of Kuranishi space of local semiuniversal deformation of a pair, we have the following theorem due to Kawamata [Cor. 4, \cite{Kaw}].

\begin{Theorem}\label{thm:deform}
For each  pair $(X, H)$  of a compact complex manifold $X$ 
and a normal crossing divisor $H$, 
there exist a germ of a complex variety
$(B, o)$ and the semiuniversal deformation of $(X, H)$  
$$
\begin{array}{ccr}
 {\cal X}  & {\hookleftarrow} & {\cal H}
 = \sum_{i=1}^r {\cal H}_i  \\
\pi \downarrow \quad  &    \swarrow \varphi &   \\
 B &  &  
\end{array}.
$$
Moreover if 
$$
H^2(X, \Theta_X(- \log H)) = \{ 0 \}, 
$$
the germ $(B, 0)$ is smooth and  the Kodaira--Spencer 
map induces an isomorphism
$$
T_0(B) \stackrel{\simeq}{\lra} H^1(X, \Theta_X(-\log H)).
$$
\end{Theorem}

The following Lemma is well known and easy to verify.  

\begin{Lemma} Let $(X, H= \sum_{i=1}^l H_i)$ be as above, 
and let $\nu:\tilde{H} = \coprod_{i=1}^r H_i \lra H$ be 
the normalization map.  
Then we have exact sequences of sheaves:
\begin{equation}\label{eq:omega}
 0 \lra \Omega^1_X \lra \Omega^1_X(\log H) \stackrel{P.R.}{\lra} 
 \nu_{*}(\oplus_{i=1}^l {\cal O}_{H_i}) \lra 0
 \end{equation}
\begin{equation}\label{eq:theta}
0 \lra \Theta_X (- \log H) \lra \Theta_X \lra \nu_*(\oplus_{i=1}^l N_{H_i/X} )
\lra 0
\end{equation}
Here the map $P.R.:\Omega^1_X(\log H) \lra 
 \nu_*(\oplus_{i=1}^l {\cal O}_{H_i}) $ is induced by the Poincar\'e 
 residue and $N_{H_i/X} = {\cal O}_X(H_i)/{\cal O}_X$ denotes 
 the normal bundle of the divisor $H_i \subset X $.
\end{Lemma}

\vspace{1cm}
\noindent
{\bf Deformation of generalized rational Okamoto--Painlev\'e pairs}

\vspace{0.5cm}

Let $(S, Y)$ be a generalized rational  Okamoto--Painlev\'{e} pair.  Recall that 
 $Y = \sum_{i=1}^r m_i Y_i $ is the anti-canonical divisor $-K_S$. 
Moreover  we set $D_i = Y_{red} = \sum_{i=1}^r Y_i$.
From now on, we will calculate some   
cohomology groups of the pair $(S, D)$ for 
Okamoto--Painlev\'e pair $(S, Y)$  which we will use later.  

Let $i:D \hookrightarrow S$ be the natural  inclusion and 
$\nu:\tilde{D} = \coprod_{i=1}^r Y_i \ra D $ the normalization map.  Set $j = i \cdot \nu $.  
First, let us consider the following Gysin exact sequence
\begin{eqnarray}\label{eq:exact}
H^1(S, \C) \lra H^1(S-D, \C) \lra H^0(\tilde{D}, \C) \stackrel{j!}{\lra} H^2(S, \C) \lra \cdots.   
\end{eqnarray}

The following lemma is important.  

\begin{Lemma}\label{lem:gysin}
 Under the same notation as above, the 
Gysin map gives an injective homomorphism 
$$
H^0(\tilde{D}, \C) \hookrightarrow H^2(S, \C).
$$
\end{Lemma}

{\it Proof.}  Since 
$$
H^0(\tilde{D}, \C) = \oplus_{i=1}^r H^0(Y_i, \C) = \oplus_{i=1}^r \C \cdot {\bf 1}_{Y_i},  
$$
and the image of the Gysin map of ${\bf 1}_{Y_i} $ is just 
the divisor class  $ Y_i  \in H^2(S, \C)\simeq \Pic(S) 
\otimes \C $, we only have to show that the classes $\{ Y_i  \}_{i=1}^r $ is lineally independent in $H^2(S, \C)$.  Looking at the intersection matrix of $\{ Y_i \}_{i=1}^r$, we easily see that only possible linear relation is 
$$
 Y = \sum_{i=1}^r m_i Y_i = {\bf 0}.
$$
On the other hand, $K_S = -Y$ and $S$ has at least one $(-1)$-smooth rational curve $E$. By adjunction formula, we see that $Y \cdot E = -(K_S) \cdot E = 1$, hence we see that $\{Y_i \}_{i=1}^r$ is linearly independent.  \qed

As a corollary to Lemma \ref{lem:gysin}, we obtain:

\begin{Corollary}\label{cor:gysin}
For a generalized rational  Okamoto--Painlev\'e pair $(S, Y)$, we have  the following.
\begin{enumerate}
\item $H^1(S -D, \C) =  0 $.

\item $H^0(S, \Omega^1_S(\log D)) =  0$.

\item $H^2(S, \Theta_S(-\log D)) =  0$.  

\item $H^2(S, \Theta_S) = 0$.
\end{enumerate}
\end{Corollary}

\noindent
{\it Proof.}  
Since $S$ is a rational surface, we have $H^1(S, \C) = 0$.
From the exact sequence (\ref{eq:exact}) and 
Lemma \ref{lem:gysin}, we have the first assertion.  Then 
from the mixed Hodge theory, we have an inclusion 
$$
H^0(S, \Omega^1_S(\log D)) \hookrightarrow H^1(S -D, \C).
$$
This proves the second assertion.  
For the third assertion, let us 
consider the Serre duality 
$$
H^2(S, \Theta_S(-\log D))^{\vee} \simeq H^0(S, \Omega^1_S(\log D) \otimes K_S).
$$
Since $K_S = {\cal O}_S(-Y)$, we have an  inclusion 
$$
H^0(S, \Omega^1_S(\log D) \otimes K_S) \hookrightarrow H^0(S, \Omega^1_S(\log D)) = \{ 0 \},
$$
This shows the third assertion. 
From the exact sequence (\ref{eq:theta}), 
we obtain the exact sequence 
$$
 \lra 
 H^2(S, \Theta_S(-\log D)) 
 \lra 
 H^2(S, \Theta_S) 
 \lra 
 H^2( D, N_{D_i/S})
 \lra 0. 
$$
Since $ \dim D = 1 $,  $ H^2( D, N_{D_i/S}) = 0$, 
hence, from  the third assertion we obtain the fourth assertion.  

\vspace{0.3cm}
The following geometric facts are very important for our purpose.  (cf.  [Lemma 3, \cite{AL}], \cite{SU}). 

\begin{Proposition}\label{prop:autom}
Let $(S,Y)$ be  a generalized rational  Okamoto--Painlev\'e pair  such that $Y$ is a divisor with normal crossing and $(S, Y)$ is not of fibered type.  
\begin{enumerate}
\item $H^0(S-D, {\cal O}^{alg} ) \simeq \C$
\item $H^0(S-D, \Theta^{alg}_{S-D}) \simeq 0 $.  Here $\Theta^{alg}_{S-D}$ 
denotes the sheaf of germs of algebraic 
regular infinitesimal automorphisms.  
\item $H^0(S, \Theta_S( H)) = 0$ for any effective divisor $H$ supported 
on $D$.  
\item $H^0(S, \Theta_S(-\log D)(H)) = 0 $ for any effective divisor $H$ supported on $D$.  
\end{enumerate}
\end{Proposition} 

{\it Proof.}  
The first assertion follows from  Proposition \ref{prop:reg}.  Since on $S$ there exists a non-zero 
rational $2$-forms $\omega_S$  which is non-degenerate on $S-D$, $\omega_S$ induces an 
isomorphism $\Theta_{S-D} \simeq \Omega^1_{S-D} $.  Hence it suffices to show that $H^0(S-D, 
\Omega^{1, alg}_{S-D})= (0) $.  
Taking a section $ \eta \in H^0(S-D, 
\Omega^{1, alg}_{S-D}) $, 
we see that $ d \eta/\omega_{S} $ is 
a regular holomorphic function on $S-D$, hence constant $c$ (cf. Propositin \ref{prop:reg}).  So this implies that $d \eta = c \cdot \omega_{S}$. On the other hand, we can easily see that $\omega_{S}$ is 
non-zero 
element in $H^2_{DR}(S-D, \C) $, hence $d \eta = 0$.  
Hence it lies in $H^0(S-D, d {\cal O}_{S-D}^{alg})$.  
Since we have an isomorphism (cf. [3.1.7.1, \cite{D}]) 
$$
H^1(S-D, \C_{S-D}) \simeq {\bf H}^1_{DR}(S-D), 
$$
and $H^1(S-D, \C_{S-D}) = (0)$ from Corollary \ref{cor:gysin}, (1), we see that $\eta$ can be written as $d f$ for 
some $f \in H^0(S-D, {\cal O}_{S-D}^{alg})$.  However since $f$ is constant (Proposition \ref{prop:reg}), we see that $\eta = d f = 0$. The last two 
assertions easily follow from the second assertion.   
\qed

The following proposition  shows that the Kuranishi space of a generalized rational Okamoto--Painlev\'e pair $(S, Y)$ is smooth  and has dimension $10 -r$ where $r$ denotes the number 
of irreducible components of $Y$.   

\begin{Proposition}
Let $(S, Y)$ be  a generalized rational  Okamoto--Painlev\'e pair  such that $D = Y_{red}$ is a simple normal crossing 
divisor and $Y \not= \tilde{A}_0$-type. Then we  have
\begin{eqnarray}
c_2(S) = \mbox{topological Euler characteristic} = 12, 
\end{eqnarray}
\begin{eqnarray}
b_2(S) = \rank H^2(S, \Z) = 10, 
\end{eqnarray}
\begin{eqnarray}\label{eq:dim}
\dim H^1(S, \Theta_S ) = 10, 
\end{eqnarray}
and
\begin{eqnarray}\label{eq:logdim}
\dim H^1(S, \Theta_S(- \log D))  = 10 - r 
\end{eqnarray}
where $r $ is  the number of irreducible components of $Y$.  
Moreover, the Kuranishi space of the local deformation of 
the pair $(S, D)$ is smooth and of dimension $10 - r$. 
\end{Proposition}
\vskp

{\it Proof.}
First, from Noether's formula, we obtain
\begin{equation}
\chi(S, \cO_S) = \frac{1}{12}( (K_S)^2 + c_2(S)).
\end{equation}
From the definition of a generalized rational  Okamoto--Painlev\'{e} pair $(S, Y)$, 
we have $K_S = -Y$ and $(K_S)^2 = Y^2 = 0$. 
Since $S$ is rational,  we have $\chi(S, \cO_S) = 1$.  Hence, we have 
$$
c_2(S) =  12.\footnote{Note that this also follows from the fact that $S$ is a 9-points blown-up of $\BP^2$.}
$$

Since $B_1(S) = 0$, the above equality implies that $B_2(S) = 12 - 2 = 10$.

 From Riemann-Roch-Hirzeburuch formula,  we obtain 
$$
\chi(S, \Theta_S) = \sum_{i=1}^2 (-1)^i \dim H^i(S,\Theta_S) 
= \frac{1}{6} (7 \cdot (K_S)^2 - 5 c_2(S)).
$$
Then again from $(K_S)^2 = 0$ we obtain 
$$
\chi(S, \Theta_S) = -\frac{5}{6} c_2(S) = - \frac{5}{6} \times 12 = -10.
$$
Moreover from Corollary \ref{cor:gysin}, (4) and Proposition \ref{prop:autom} 
we obtain $H^i(S, \Theta_S) = 0$ for $i = 0, 2$, and hence
\begin{equation}\label{eq:theta10}
\dim H^1(S, \Theta_S) = 10.  
\end{equation}
For $H^1(S, \Theta_S(- \log D))$, consider the exact sequence 
\begin{equation}\label{eq:theta2}
0 \lra \Theta_S (- \log D) \lra \Theta_S \lra \nu_*( \oplus_{i=1}^r N_{Y_i/S} )
\lra 0.
\end{equation}
Remember $D = Y_{red} = \sum_{i=1}^r Y_i$. Since $K_S \cdot Y_i = 0$ and 
$Y_i \simeq \BP^1$ by assumption, 
the adjunction formula shows that
$$
\deg N_{Y_i/S} = -2  \ \   \mbox{or} \  \  N_{Y_i/S} \simeq \cO_{\BP^1}(-2).
$$
Then since $H^0(\BP^1, \cO_{\BP^1}(-2)) = 0$, $H^1(\BP^1, \cO_{\BP^1}(-2)) 
\simeq \C$ and $H^2(S, \Theta_S(-\log D)) = 0$, 
we have the following exact sequence 
\begin{equation}
0 \lra H^1(S, \Theta_S(-\log D)) \lra H^1(S, \Theta_S) \lra \oplus_{i=1}^r \C[Y_i] \lra  0 
\end{equation}
This  implies that the assertion (\ref{eq:logdim}) holds.
The last assertion follows from Theorem \ref{thm:deform} and the fact that 
$H^2(S, \Theta_S( - \log D)) = 0$.

\vskp
\begin{table}\label{tab:def}
\begin{center}
{\bf Table of the deformation of generalized rational  Okamoto--Painlev\'e pairs.}

\vspace{0.5cm}
\begin{tabular}{||c||c|c|c|c|c|c|c|c|c|c||} \hline
    &   & & & & & & & &  \\
$Y$ & $\tilde{E_8}$ & $\tilde{D_8}$  & $ \tilde{E_7}$ & 
$ \tilde{D_7}$ & 
$ \tilde{D_6}$ & $\tilde{E_6}$ & $\tilde{D_5}$ & $\tilde{D_4}$ & $\tilde{A}_{r-1}$, $ r \geq 2 $ \\
    &  & & & & & & & & \\ \hline
        & &   & & & & & & &  \\
Number of components of $Y$& $9 $ & $9$ & $ 8 $ & $8$ &  $ 7 $ & $ 7$  
 & $ 6 $ & $5$ &  $r$ \\ 
    &  & & & & & & &  &  \\ \hline
   &  & & & & & & &  & \\ 
$\dim H^1(S, \Theta_S(-\log D))$ & $1$ & $1$ & $2$ & $2$ & $3$ & $3$ 
& $4$  &  $5$ & $10 - r$\\
    & & & & & & & & & \\ \hline
    & & & & & & & & & \\ 
Painlev\'{e} equation & $P_{I}$ & $P^{\tilde{D}_8}_{III}$ & $P_{II}$ & $P^{\tilde{D}_7}_{III}$ & $P^{\tilde{D}_6}_{III}$ & $P_{IV}$ 
& $P_{V}$  &  $P_{VI}$ & none\\  
   &   & & & & & & & &  \\ \hline
\end{tabular}
\vspace{0.3cm}
\caption{}
\end{center}
\end{table}

\vspace{1cm}
\section{Local cohomology sequences and Time variables}
\label{sec:local}

Let $(S, Y)$ be a generalized rational  Okamoto--Painlev\'{e} pair and set $D = Y_{red}$. Moreover, in this section, 
we assume that
\begin{enumerate}

\item $(S, Y)$ is  of non-fibered type and 

\item  $Y_{red}$ is a normal crossing divisor with 
at least two irreducible components,  i.e. $(r \geq 2)$ 
so that all irreducible components of $Y_{red}$ are 
smooth rational curves.  
\end{enumerate}

In what follows, $\cO_S$ and  $\cO_{S-D}$ denote the sheaves
of germs of algebraic regular functions on $S$ and $S - D$ 
respectively. Moreover all sheaves of $\cO_S$-modules are 
considered in algebraic category unless otherwise stated.   
Let us consider the following exact sequence of local 
cohomology groups ([Corollary 1.9, [Gr]])
\begin{eqnarray}
H^0(S, \Theta_S(- \log D )) &  \ra H^0(S-D, \Theta_{S}( - \log D)) & \ra H^1_D(\Theta_S( - \log D)) \ra \\ 
 H^1(S, \Theta_S( - \log D)) & 
 \stackrel{\res}{\ra} H^1(S - D, \Theta_S( - \log D))  & .
\end{eqnarray}

Since $(S, Y)$ is of non-fibered type, from (2), Proposition  \ref{prop:autom}, we see that  
$$
H^0(S-D, \Theta_S (- \log D)) =  H^0(S -D, \Theta_S) = \{ 0 \}.
$$ 

Hence, we have the following 

\begin{Proposition}\label{prop:local}
For a generalized rational Okamoto--Painlev\'e pair of non-fibered type, we have the following exact sequence: 
\begin{equation}\label{eq:local}
\begin{array}{cccccc}
0 \ra 
&  H^1_D(\Theta_S( - \log D)) 
&  \ra 
&  H^1(S, \Theta_S( - \log D))
 & \stackrel{\res}{\ra} 
 & H^1(S - D, \Theta_S( - \log D)) \\
\end{array}
\end{equation}
\end{Proposition}

The following theorem is proved in \cite{T}. 

\begin{Theorem}\label{thm:time}
 Let $(S, Y)$ be  a generalized rational Okamoto-Painlev\'e pair $(S, Y)$ with the condition above. Moreover $D = Y_{red}$ is of additive type. Then  we have 
\begin{equation}
\dim H^0 (D, \Theta_S(- \log D) \otimes N_{D} ) = 1.
\end{equation}
Here  we put $N_D = {\cal O}_S(D)/{\cal O}_S$.  

Since we have 
a natural inclusion 
$$
H^0 (D, \Theta_S(- \log D) \otimes N_{D} ) \hookrightarrow H^1_D(\Theta_S(-\log D) ), 
$$ 
we obtain 
\begin{equation}\label{eq:time}
\dim H^1_D(\Theta_S(-\log D) ) \geq  1.
\end{equation}
\end{Theorem}

This theorem plays an important role to understand the 
Painlev\'e equation related to $(S, Y)$.  
We will not investigate the further 
structure of local cohomology here. Instead,   we propose 
the following
\begin{Conjecture}\label{conj:local}
Under the same notation and assumption as in Theorem \ref{thm:time},  
\begin{equation}
H^1_D(\Theta_S(-\log D) ) \simeq H^0(D, \Theta_S(- \log D) \otimes N_{D} ) \simeq \C.
\end{equation}
\end{Conjecture}

From the exact sequence (\ref{eq:local}), we see that 
the subspace $ H^1_D(S, \Theta_S(-\log Y))$ of 
$H^1(S, \Theta_S(- \log Y))$ coincides with the kernel of $\mu$. This implies that:
\begin{equation*}
H^1_D(S, \Theta_S(- \log D)) \simeq 
\left\{ \mbox{\begin{tabular}{l}
Infinitesimal deformations of $(S, D)$ whose 
restriction  \\
to $S-D$ induces the  trivial deformation
\end{tabular}} \right\}.
\end{equation*}

In \S \ref{sec:globaltoham},  
we will show that  any direction 
corresponding
 to a non-zero element of the local cohomology 
group $H^1_D(S, \Theta_S(- \log D))$ induces a 
differential equation (at least locally) by using 
\v{C}ech coboundaries.  

At this moment, we can not  prove Conjecture \ref{conj:local} in the full generality. 
However,  we see that 
the one dimensional 
vector subspace $ H^1(D, \Theta_S(-\log D)\otimes N_D ) $ of $H^1_D(\Theta_S(-\log D) ) \subset H^1(\Theta_S(-\log D) )$ 
really corresponds to the time variable $t$ in the known 
Painlev\'e equation.  
It is unlikely that we will have more time 
variables, so this gives an  evidence of Conjecture \ref{conj:local}.  

Let us explain the strategy of proving Theorem \ref{thm:time} in \cite{T}. Recall that 
$$
H^1_D(S, \Theta_S( - \log D)) = \varinjlim {\rm Ext}^1({\cO_{nD}}, \Theta_S( - \log D)) 
$$
where $\cO_{nD} = \cO_S/\cO_S(-n D)$ (cf. [Theorem 2.8, [Gr]]). 

On the other hand, since $\Theta_S( - \log D)$ is a 
locally free $\cO_S$-module, we see that
\begin{equation}
{\cal H}om({\cO_{nD}}, \Theta_S( - \log D)) = 0, 
\end{equation}
and 
\begin{equation}
{\cal E}xt^1({\cO_{nD}}, \Theta_S( - \log D)) =  \Theta_S( - \log D) \otimes N_{nD}, 
\end{equation}
where $N_{nD} = \cO_S(nD)/\cO_S$.  
By an argument using a spectral sequence, we see that 
\begin{equation}
H^1_D(S, \Theta_S( - \log D)) = \varinjlim H^0( \Theta_S( - \log D) \otimes N_{nD}) 
\end{equation}
Hence, we have a natural inclusion
\begin{equation}
H^0( \Theta_S( - \log D) \otimes N_{D}) \hookrightarrow H^1_D(S, \Theta_S( - \log D)).  
\end{equation}

\begin{Lemma} Let $(S, Y)$ be a generalized rational Okamoto--Painlev\'e pair as above and set $D= Y_{red}$.  
Then we have the following exact sequences
\begin{equation}\label{eq:thetanormal}
0 \lra \Theta_D \otimes N_D \lra \Theta_S \otimes N_D \lra 
\nu_*(\oplus_{i=1}^r N_{Y_i/S})\otimes N_D \lra 0.
\end{equation}
\begin{equation}\label{eq:thetalog}
0 \lra \nu_*(\oplus_{i=1}^r N_{Y_i/S}) 
 \lra \Theta_S(- \log D) \otimes N_D \lra  \Theta_D \otimes N_D
 \lra 0. 
\end{equation}
Here $\Theta_D$ denotes the tangent sheaf of $D$ and $\nu:\tilde{D} \lra D$ the normalization map. 
\end{Lemma}

{\it Proof.}  The first exact sequence 
(\ref{eq:thetanormal}) follows from [(1.9), \cite{B-W}].  

Let us consider the following  diagram:  
$$
\begin{array}{ccccccc}
&0&&0&& {\cal \ker \lambda}& \\
 & \downarrow & & \downarrow  & & \downarrow & \\
 &&&&&& \\
0 \lra & 
\Theta_S(- \log D) & 
\lra & 
\Theta_S(- \log D) \otimes \cO_S(D) & 
\lra & 
\Theta_S(- \log D) \otimes N_D & 
\lra 0 \\
 & \downarrow & & \downarrow  & & \quad \downarrow \lambda & \\
0 \lra & 
\Theta_S& 
\lra & 
\Theta_S \otimes \cO_S(D) & 
\lra & 
\Theta_S \otimes N_D & 
\lra 0 \\
 & \downarrow & & \downarrow  & & \downarrow & \\
 & \nu_*(\oplus_{i=1}^r N_{Y_i/S}) & \stackrel{\mu}{\lra} & 
 \nu_*(\oplus_{i=1}^r N_{Y_i/S}) \otimes N_D  & \lra & \coker \lambda  
  &  \lra 0 \\
   & \downarrow & & \downarrow  & & \downarrow  & \\
&0&&0&& 0 &. \\
\end{array}
$$
By the snake lemma, we obtain the exact sequence 
$$
0 \lra \ker \lambda \lra \nu_*(\oplus_{i=1}^r N_{Y_i/S}) \stackrel{\mu}{\lra} \nu_*( \oplus_{i=1}^r N_{Y_i/S}) \otimes N_D   \lra  \coker \lambda  
    \lra 0.     
$$
From a local consideration, we see that the map $\mu$ is the  zero map, hence
$$
\ker \lambda \simeq \nu_*(\oplus_{i=1}^r N_{Y_i/S}),  \quad 
\nu_*(\oplus_{i=1}^r N_{Y_i/S}) \otimes N_D   \simeq  \coker \lambda.
$$

Moreover since $\im \lambda \simeq \ker [ \Theta_S \otimes N_D \lra  \nu_*( \oplus_{i=1}^r N_{Y_i/S}) \otimes N_D ]$, from the exact sequence (\ref{eq:thetanormal}), we obtain the exact sequence  (\ref{eq:thetalog}).

\qed.

From the exact sequence (\ref{eq:thetalog}), one can obtain 
\begin{equation}
H^0(\oplus_{i=1}^r N_{Y_i/S}) \lra  
 H^0(\Theta_S(- \log D) \otimes N_D)  \lra H^0( \Theta_D \otimes N_D) \stackrel{\delta}{\lra} H^1(\oplus_{i=1}^r N_{Y_i/S}).
\end{equation}
where $\delta$ denotes the connected homomorphism.  

Note that since $N_{Y_i/S} = \cO_{Y_i}(-2)$, we have 
$$
H^0(\oplus_{i=1}^r N_{Y_i/S}) = \{0 \}, \quad  H^1(\oplus_{i=1}^r N_{Y_i/S} ) \simeq \C^r. 
$$

Moreover,  one can easily see that
$$
 \Theta_D \simeq \nu_*(\oplus \Theta_{Y_i} ( -t_i))  \simeq 
 \nu_*(\oplus \cO_{Y_i} (2 - t_i)) 
$$
where $t_i$ is the number of intersections of $Y_i$ with the other 
$Y_j$s.  On the other hand, since $D \cdot Y_i = t_i - 2$, we see that
$$
H^0(\Theta_D \otimes N_D)  \simeq H^0(\oplus_{i=1}^r \cO_{Y_i}), 
$$
hence 
\begin{equation}\label{eq:inclusion}
 H^0(\Theta_D \otimes N_D)
 \simeq \C^r. 
\end{equation}

From this, the connecting  homomorphism $\delta$  
\begin{equation}\label{eq:delta}
\delta : H^0(\Theta_D \otimes N_D)  \lra   \oplus_{i=1}^r H^1(N_{Y_i/S}) 
\end{equation}
can be identified with  a linear map  $\delta:\C^r \lra \C^r$ and 
we have an isomorphism 
\begin{equation}
 H^0(D, \Theta_S(- \log D) \otimes N_D ) \simeq \ker \delta. 
\end{equation}

The following proposition is the main theorem of \cite{T}.  

\begin{Proposition}\label{prop:main}
 Let $(S, Y)$ be as in Theorem \ref{thm:time}.  
A matrix representation of the linear map  $\delta:H^0(\Theta_D \otimes N_D)  \lra   \oplus_{i=1}^r H^1(N_{Y_i/S})$  in (\ref{eq:delta}) is 
equal to the $\pm$ of the intersection matrix of $D = \sum_{i=1}^r Y_i $, that is, 
$$
\delta = ( (Y_i \cdot Y_j))_{1\leq i, j \leq r}.
$$
Since the intersection matrix $( (Y_i \cdot Y_j))_{1\leq i, j \leq r}$  has exactly one-dimensional kernel 
corresponding to the space of $ Y = \sum_{1=1}^r m_i Y_i$, we have 
$$
H^0(D, \Theta_S(- \log D)\otimes N_D ) \simeq  \ker \delta = \C.  
$$ 
\end{Proposition}

\section{Reviews on Kodaira--Spencer theory}\label{sec:KS}

In this section, we review on Kodaira--Spencer theory of 
complex analytic deformation.  A main reference is 
\cite{KS}.  

Let $X$ be a compact complex manifold  of dimension $n$.  
We can take a locally finite open covering $\{ U_i \}_{i \in I}$ of $X$ such that each open subset $U_i$ admits  local coordinates $ \z_i = (z_i^{1}, \cdots, z_i^{n})$: 
$$
X = \cup_{i \in I} U_i.
$$
For a point  $ p \in U_i \cap U_j$,  we have two local 
coordinates $\z_i(p)$ and $\z_j(p)$ whose coordinate transformation are given by 
$$
\z_i =(z_i^{1}, z_i^{2}, \cdots, z_i^n) = {\bf f}_{ij}(\z_j), 
$$
or more precisely for $\alpha =1, \cdots, n $, 
\begin{equation}
  z_i^{\alpha} = f_{ij}^{\alpha}(z_j^1, \cdots. z_j^n).  
\end{equation}
Here $f_{ij}^{\alpha}(\z_j)$ are holomorphic functions defined on 
$ U_i \cap U_j $.  
Note that one can give  the compatibility conditions for $ U_i \cap U_j \cap U_k \not= \emptyset $
\begin{equation}
f_{ik}^{\alpha}(\z_k) = f_{ij}^{\alpha}(f_{jk}^{1}(\z_k), 
\cdots, f_{jk}^n(\z_k))
\end{equation}

Complex structure of $X$ can be deformed 
by changing these coordinate transformations keeping the 
compatibility conditions.  

Let $\cB$ be a complex manifold  with a (global) coordinate system $(t_1, \cdots, t_m) $ and a specific marked point $\0 = (0, \cdots, 0)  \in \cB$.  We may think that $\cB$ is an 
affine variety or a complex analytic small open  ball around the origin.  
\begin{Definition}{\rm 
A deformation of $X$ with a parameter space $\cB \ni(t_1, \cdots, t_m)  $ is a proper smooth  holomorphic map 
$\pi:{\cal X} \lra \cB$ such that the following 
diagram is commutative:
$$
\begin{array}{cccc}
 {\cal X}  & \stackrel{\iota}{\hookleftarrow} & X_0 & \simeq X
       \\
\pi \downarrow \quad  &    &    \downarrow &  \\
 \cB &  \ni  &  \0  & .
\end{array}
$$
}
\end{Definition} 

\begin{Definition}{\rm 
A deformation $\pi:{\cal X} \lra \cB$ of $X$ 
 is said to have {\em a 
finite covering relative to  $\cB$} if ${\cal X}$ is covered by $\{ \tilde{U_i} =  U_i \times \cB  \}$ 
such that the following diagram is commutative:
$$
\begin{array}{ccc}
 {\cal X} & = &  \cup_{i \in I} (U_i \times \cB)  \\
\pi \downarrow \quad  &   &     \downarrow    \\
 \cB &  =  & \cB  \\. 
\end{array}
$$
}
\end{Definition} 

Let us assume that $\pi:{\cal X} \lra \cB$ 
has a finite covering relative to $\cB$ and take 
the local coordinate of $\tilde{U_i} = U_i \times \cB$ 
by $(z_i^{1}, \cdots, z_i^n, t_1, \cdots, t_m ) $.  
The coordinate transformation for $\tilde{U_i} \cap \tilde{U_j}$ is given by 
$$
z_i^{\alpha} = f_{ij}^{\alpha}(z_j^1, \cdots, z_j^n, t_1, \cdots, t_m ). 
$$
We may assume that for ${\bold t} = \0$, we have $f_{ij}^{\alpha}(z_j^1, \cdots, z_j^n, 0, \cdots, 0) = f_{ij}^{\alpha}(z_j^1, \cdots, z_j^n)$. 

Now we can introduce the {\em Kodaira--Spencer class} 
of the deformation  $\pi:{\cal X} \lra \cB$  for 
each ${\bold t} \in \cB$.

For simplicity we assume that $\cB$ is one dimensional,  
hence $t=t_1$ is the global parameter of $\cB$.  
Let $h$ be a holomorphic function on an open subset $V$ of ${\cal X}$.  
Then on $\tilde{U_i} \cap V$, $h$ is a function 
in a local coordinate $h(z_i^1, \cdots, z_i^n, t)$.   
Assume that $\tilde{U_i}\cap \tilde{U_j} \cap V \not= 
\emptyset$.  Regarding as
$$
h(\z_i, t) = h(f_{ij}^1(\z_j, t), \cdots, f_{ij}^n(\z_j, t), 
t), 
$$
from the chain rule, we obtain   
\begin{equation}
\left( \frac{\partial h}{\partial t}\right)_{j} = 
\left(\frac{\partial h}{\partial t}\right)_{i} + \sum_{\alpha=1}^n 
\frac{\partial f_{ij}^{\alpha}(\z_j, t) }{\partial t} \frac{\partial h}{\partial z_i^{\alpha}}.  
\end{equation}

This implies that,  as a vector field on $\tilde{U_i} \cap 
\tilde{U_j}$,  we have the following identity:
\begin{equation}\label{eq:lift}
\left( \frac{\partial }{\partial t}\right)_{j} = 
\left(\frac{\partial }{\partial t}\right)_{i} + \sum_{\alpha=1}^n 
\frac{\partial f_{ij}^{\alpha}(\z_j, t) }{\partial t} \frac{\partial }{\partial z_i^{\alpha}}.  
\end{equation}

Let us set $ \{\theta_{ij}(t) \} $ by 
\begin{equation}\label{eq:KSclass}
\theta_{ij}(t) =  \sum_{\alpha=1}^n 
\frac{\partial f_{ij}^{\alpha}(\z_j, t) }{\partial t} \frac{\partial }{\partial z_i^{\alpha}}.  
\end{equation} 
From the compatibility conditions for $ \tilde{U_i} 
\cap \tilde{U_j} \cap \tilde{U_k} \not= \emptyset $
\begin{equation}
f_{ik}^{\alpha}(\z_k, t) = f_{ij}^{\alpha}(f_{jk}^{1}(\z_k, t), 
\cdots, f_{jk}^n(\z_k, t), t), 
\end{equation} 
we obtain the identity 
$$
\theta_{ik}(t) = \theta_{ij}(t) + \theta_{jk}(t).
$$
This implies that 
$
\{\theta_{ij}(t) \}
$ defines a 
\v{C}ech 1-cocycle with values in $\Theta_{X_t}$, hence defines a cohomology class 
$$
\theta(t) \in H^1(X_t, \Theta_{X_t}),  
$$  
which is called the {\em Kodaira--Spencer class} .  

If the dimension of $\cB$ is greater than one, we can define
the cohomology class for each $\frac{\partial}{\partial t_{\mu}}$.  
More precisely one can  define  the {\em Kodaira-Spencer map}
\begin{eqnarray}
\rho:T_{{\bold t}}(\cB) & \lra & H^1(X_t, \Theta_{X_t}) \\
      v  & \mapsto &  \rho(v) = \theta_v({\bold t}) 
\end{eqnarray}
by 
$$
\theta_{v, ij} = \{ \theta_{v, ij}({\bold t}) = \sum_{\alpha=1}^{n} v( f_{ij}^{\alpha}(\z_j, {\bold t})) \frac{\partial}{\partial z_i^{\alpha}} \} . 
$$
Here for 
$$
v = \sum_{\mu=1}^{m} A_{\mu}(\bt) \frac{\partial}{\partial t_{\mu}}, 
$$ we set 
$$
v( f_{ij}^{\alpha}(\z_j, {\bold t})) = \sum_{\mu=1}^{m} A_{\mu} \frac{\partial f_{ij}^{\alpha}(\z_j, {\bold t}) }{\partial t_{\mu}}.
$$

\begin{Definition} 
{\rm A deformation $\pi:{\cal X} \lra \cB$ 
is called {\em locally trivial}, 
if for each point $\bt \in \cB$ 
there exists an  open neighborhood $I$ of $\bt$ such that 
${\cal X}_{|I} \lra I $ is complex analytically isomorphic to the product ${\cal X}_{\bt} \times I $.
}
\end{Definition}

\begin{Proposition} $($\cite{KS}$)$ 
Let $ \pi:{\cal X} \lra \cB$ be a deformation of 
a compact complex manifold with parameter space $\cB \ni 
\bt=(t_1, \cdots, t_m)$.  If for every point $\bt \in \cB$ 
$\dim H^1({\cal X}_{\bt}, \Theta_{ {\cal X}_{\bt}})$ is constant and 
 Kodaira--Spencer map $\rho$ is the zero map, then $ \pi:
 {\cal X} \lra \cB $ is a locally trivial fibration.  
\end{Proposition}

For a  proof in detail, 
we refer the reader to [Theorem 5.1, \cite{KS}].  Since we will use the idea of the proof later, we explain an outline  of the proof of  theorem when 
$\dim \cB = 1$. By replacing  $\cB$  with a neighborhood  of $\bt \in \cB$, we may  assume that 
a deformation $\pi:{\cal X} \lra \cB $ has  a finite covering $\{ \tilde{U}_i = U_i \times \cB \}$ relative to $\cB$.  
Then the  Kodaira--Spencer class 
$$
\rho(\frac{\partial}{\partial t}) = \theta(t) \in H^1({\cal X}_t, \Theta_{{\cal X}_t}). 
$$ 
is represented by \v{C}ech cocycles  $\{ \theta_{ij}(t) \}$ given in (\ref{eq:KSclass}).  Since $\theta(t)$ is cohomologus to zero, 
for each $t$ we can find 
$$
\theta_i(t) \in \Gamma(\tilde{U_i}\cap {\cal X}_t, \Theta_{U_i}). 
$$
such that
$$
\theta_{ij}(t) = \theta_j(t) - \theta_i (t) \quad \mbox{on} 
\quad 
\tilde{U}_i \cap \tilde{U}_j \cap {\cal X}_t
$$
From (\ref{eq:lift}), we obtain the following identities  of vector fields  on each $\tilde{U}_i \cap \tilde{U}_j \cap {\cal X}_t $
\begin{equation}
\left( \frac{\partial }{\partial t}\right)_{j} = 
\left(\frac{\partial }{\partial t}\right)_{i} + (\theta_j(t)
- \theta_i (t) ),   
\end{equation}
and hence 
\begin{equation}\label{eq:vf}
\left( \frac{\partial }{\partial t}\right)_{j} - \theta_j(t)  = 
\left(\frac{\partial }{\partial t}\right)_{i} - \theta_i(t). 
\end{equation}

At this moment, it is not obvious that 
the dependence of 
\begin{equation}
\theta_{i}(t) = \sum_{\alpha=1}^{n} \theta_{i}^{\alpha}(\z_i, t) 
\frac{\partial }{\partial z_i^{\alpha}}
\end{equation}
with respect to $t$ 
is in $C^{\infty}$ class.  However 
under the condition that $\dim H^1({\cal X}_t, \Theta_{{\cal X}_t}) $ is constant on $\cB$,  one can prove that 
$ \theta_i(t) $ can be chosen as a
vector field on $\tilde{U}_i= U_i \times \cB $ in
$C^{\infty}$ class. 

 This implies that the vector field 
 \begin{equation}
  \{ \left( \frac{\partial }{\partial t}\right)_{i} - \theta_i(t) \}_{i \in I}
  \end{equation}
   on $ \tilde{U}_i $ can be glued together and defines 
a global $C^{\infty}$--vector field, say, $ \tilde{v} $ on the total space 
${\cal X}$. 
We see that $\tilde{v}$ is a lift of vector field 
$\frac{\partial }{\partial t }$ by $\pi$.  Then on each 
open set $\tilde{U}_i$, we can consider the ordinary 
differential equation 
\begin{equation} \label{eq:def}
 \frac{d z_i^{\alpha}}{d t}  = - \theta_i^{\alpha}(\z_i, t) 
 \quad \alpha = 1, \cdots, n. 
\end{equation}
And these set of differential equations can be patched together on whole ${ \cal X}$.  Starting from an initial conditions $(a_1, \cdots, a_n, t_0) \in {\cal X}_{t_0}$, the solution $(z_1(a_j, t), \cdots, z_n(a_j, t))$ of differential 
equation (\ref{eq:def})  defines a 
$C^{\infty}$--curve which 
is transversal to each fiber ${\cal X}_t$.  Then the 
whole solutions of  (\ref{eq:def}) define a foliation on ${\cal X}$ and   define  
$C^{\infty}$-defeomorphisms  $ \varphi_t: {\cal X}_{0} \stackrel{\simeq}{\ra} {\cal X}_t $. Moreover, one can show that 
this defeomorphism $\varphi_t$  is a complex biholomorphic 
morphism for each $t \in \cB$.  

This implies  the following. If we have a family of compact complex manifolds $\pi:{\cal X} \lra \cB$  with  a parameter $t \in \cB$ such that the  Kodaira--Spencer map $\rho_t$ is zero, we will obtain a differential equation as in (\ref{eq:def}) defined on the total space ${\cal X}$.   

Summarizing these, we have the following implications (cf. Figure \ref{fig:KScompact}). 

\medskip
\begin{figure}[h]
$$
\begin{array}{c}
\fbox{\begin{tabular}{l} Deformation $\pi: \cX \lra \cB$ of complex  manifolds \\ with zero Kodaira--Spencer map \end{tabular}} \\
  \\
\Downarrow \\
   \\
\fbox{There exists a vector field $ \tilde{v}$ on ${\cal X}$ which is a lift of $\frac{\partial}{\partial t}$. } \\
  \\
\Downarrow \\
  \\
\fbox{Differential Equation on $\pi:{\cal X} \lra \cB $ given by $\tilde{v}$} \\
  \\
\Downarrow \\
  \\
\fbox{Local trivializations of the deformation $\cX \lra \cB $ }.
\end{array}
$$
\caption{}
\label{fig:KScompact}
\end{figure}
\medskip

\section{Global Deformations of Okamoto--Painlev\'e pairs}
\label{sec:globaldef}

\noindent
{\bf Affine coverings and Symplectic Structures on $S - D$ }

\vspace{0.3cm}

Let $(S, Y)$ be a generalized rational Okamoto--Painlev\'e 
pair.  Then by definition, $S$ has a rational 2-form $\omega_Y$ whose pole divisor is $Y$.  Setting $D = Y_{red}$, the 
rational 2-form $\omega_Y$ induces  a non-degenerate holomorphic 2-forms on the 
open surface $S -D$, hence induces a holomorphic symplectic structure on $S -D$.    

In \cite{O1}, Okamoto introduced 
 the space of initial conditions of 
Painlev\'e equation of each type, which  can be 
written as $ S - D $ for an  Okamoto-Painlev\'e pair $(S,Y)$.  The main 
reason why Painlev\'e equations can be  written as  
Hamiltonian systems is this holomorphic symplectic structure.    
For Painlev\'e equations $P_J$, ($J = II, III, IV, V, VI$), Takano et al. \cite{ST}, \cite{MMT} constructed  a good  
family of Okamoto--Painlev\'e pairs $(S_{{\boldsymbol \alpha}, t}, Y_{{\boldsymbol \alpha}, t} ) $  depending on the time variable and   a system of auxiliary parameters 
${\boldsymbol \alpha}=(\alpha_1, \cdots, \alpha_s)$ appeared in each Painlev\'e equation. 

Summarizing  results in \cite{ST}, \cite{MMT}, let us explain the situation of spaces of initial 
conditions of classical Painlev\'e equations in the way of 
our setting. 
  Let $R=R(Y)$ be a type of the root systems corresponding  to 
  a Painlev\'e equation. 
  Then there exist an affine open subset $\cM_R$ of $\C^s = \Spec \C[\balpha]= \Spec \C[\alpha_1, \cdots, \alpha_s]$, an affine open subset $\cB_R$ 
  of $\C = \Spec \C [t]$ and the following
  deformation of non-singular pair  
\begin{equation}\label{eq:family}
\begin{array}{ccl}
 {\cal S} &  \hookleftarrow & {\cal D}      \\
\pi \downarrow \hspace{0.3cm} &    \swarrow &  \varphi  \\
 \cM_R \times \cB_R  &  & 
\end{array} 
\end{equation}
where ${\cal S} \lra \cM_R \times \cB_R$ is a smooth family of rational surface and ${\cal D} \hookrightarrow {\cal S}$ is 
a normal crossing divisor.  In order to relate this diagram to Okamoto--Painlev\'e pair, let $ \Omega^2_{{\cal S}/\cM_R \times \cB_R}(*{\cal D})$ denote 
the sheaf of germs of relative rational two forms on 
${\cal S}$ which have  poles only  along ${\cal D}$.  
There exists a section 
$$
\omega_{\cal S} \in \Gamma({\cal S}, \Omega^2_{{\cal S}/\cM_R \times \cB_R}(*{\cal D}) )
$$
which induces a rational $2$-form $\omega_{{\cal S}_{\balpha, t}}$ for each fiber ${\cal S}_{\balpha, t}$.  The pole divisor $\omega_{\cal S}$ is denoted by ${\cal Y}$, and 
with suitable choice of $\omega_{\cal S}$ 
we may assume that  each fiber $({\cal S}_{\balpha, t}, {\cal Y}_{\balpha, t})$ is an Okamoto--Painlev\'e pair of type $R = R(Y)$ and ${\cal Y}_{red} = {\cal D}$. (Note that on ${\cal S} - {\cal D} $ the relative rational 2-form $\omega_{\cal S}$ is 
holomorphic and non-degenerate on each fiber ${\cal S}_{\balpha, t} - {\cal D}_{\balpha, t}) $.  Assuming that the family (\ref{eq:family})  is effectively parameterized and semiuniversal  at each point 
of $\cM_R \times \cB_R$, or equivalently the Kodaira--Spencer map
\begin{equation}
\rho:T_{\balpha, t}(\cM_R \times \cB_R) \lra H^1({\cal S}_{\balpha, t}, \Theta_{{\cal S}_{\balpha, t}}(- \log {\cal D}_{\balpha, t}))
\end{equation}
is an isomorphism at each point, we have the equality
\begin{equation}
\dim \cM_R =  \dim H^1({\cal S}_{\balpha, t}, \Theta_{{\cal S}_{\balpha, t}}(- \log {\cal D}_{\balpha, t})) - 1 . 
\end{equation}
\vspace{0.5cm}

\begin{table}[h]
\begin{center}
{\bf The dimensions of $\cM_R$ for generalized \\ rational Okamoto--Painlev\'e pairs.}

\vspace{1cm}
\begin{tabular}{||c||c|c|c|c|c|c|c|c|c||} \hline
    &  & & &  & & & &  \\
$R= R(Y)$ & $\tilde{E_8}$ & $\tilde{D_8}$  & $\tilde{E_7} $& $\tilde{D_7}$  & 
$\tilde{D_6}$ & $\tilde{E_6}$ & $\tilde{D_5}$ & $\tilde{D_4}$ \\
    &  & & & & & & &  \\ \hline
    &  & & & & & & &  \\
Painlev\'{e} equation & $P_{I}$ &  $P^{\tilde{D}_8}_{III}$  & $P_{II}$ & $P^{\tilde{D}_7}_{III}$ & $P^{\tilde{D}_6}_{III}$ & $P_{IV}$
& $P_{V}$  &  $P_{VI}$ \\
    &   & & & & & & & \\ \hline
        &   & & & & & & & \\ 
$s = s(R) = \dim \cM_R$ &   $0$  & $0$&   $1$ & $1$ & $2$ & $2$  
& $3$   & $4 $ \\
(= $\sharp$ of auxiliary parameters.) &       &   & & & & & &  \\ \hline
\end{tabular}
\vspace{0.3cm}
\caption{}
\label{tab:aux}
\end{center}
\end{table}

(Note that the Okamoto--Painlev\'e pairs of type $\tilde{D_8}, \tilde{D_7}$ did not appear in the classical literatures
(cf. \cite{O1}, \cite{O2}, \cite{MMT}).)

More notably, Takano et al. \cite{ST}, \cite{MMT} constructed  an affine open covering $\{ \tilde{U_i} \}_{i=1}^{l + m} $ 
  of  ${\cal S}$ for the classical Okamoto--Painlev\'e pair of Painlev\'e equation $P_J$ ($J=II, \cdots, VI$),  which is relative to $\pi$ and so that 
$$
  \tilde{U_i} = \cM_R \times \cB_R \times U_i 
$$
where $U_i = \Spec \C[x_i, y_i] \simeq \C^2$.  
Moreover, we may assume that 
$\{ \tilde{U_i} \}_{i=1}^{l} $  covers 
${\cal S} -{\cal D}$ and for $ 1 \leq i \leq l$, we have 
\begin{equation}
\omega_{\cS|\tilde{U}_i} = dx_i \wedge d y_i
\end{equation} 

In this sense, the restricted morphism 
$$
\pi:{\cal S} -{\cal D} \lra \cM_R \times \cB_R
$$
is a deformation of open symplectic surfaces.

By using the results in Appendix B of \cite{Sakai}, 
we can generalize the result of Takano, et al. as follows.

\begin{Proposition} \label{prop:global}
 Let $R = R(Y)$ be one of  types of the root systems  appeared in Proposition \ref{prop:classf} which is additive type, so that 
 $$ 
 \dim H^1_D(\Theta_S(-\log D)) \geq 1
$$ 
for corresponding generalized rational Okamoto--Painlev\'e pair $(S, Y)$ (cf. Theorem \ref{thm:time}).    (That is, $R \not= \tilde{A}_{r-1}$).   
 Moreover denote by $r$ the number of irreducible components of $D = Y_{red}$.  
 
Let $ \cM_{R}$ be an  
affine open subscheme  in $\C^s = 
\Spec \C[ \alpha_1, \cdots, \alpha_s]$ of 
dimension $s =s(R) = 9 - r$  and $\cB_R $ be an affine open 
subscheme of $\C = \Spec \C[t]$.    Then there exists 
the following commutative diagram satisfying the 
conditions below.  
\begin{equation}\label{eq:globalf} 
\begin{array}{ccl}
 {\cal S} &  \hookleftarrow & {\cal D}      \\
\pi \downarrow \hspace{0.3cm} &    \swarrow &  \varphi  \\
 \cM_R \times \cB_R  &  &.  
\end{array} 
\end{equation}
\begin{enumerate}
\item The above diagram is a deformation of non-singular pair of projective surfaces and normal crossing divisors in the sense of Definition \ref{thm:deform}

\item There exists a rational relative 2-form 
$$
\omega_{\cal S} \in \Gamma({\cal S}, \Omega^2_{{\cal S}/\cM_R \times \cB_R}(*{\cal D}) ) 
$$
which has poles only along ${\cal D}$.  
\item If we denote by ${\cal Y}$ the pole divisor of $\omega_{\cal S}$, then for each point $(\balpha, t) \in \cM_R \times \cB_R$,  $({\cal S}_{\balpha, t}, {\cal Y}_{\balpha, t})$ is a generalized Okamoto--Painlev\'e pair of type $R = R(Y)$ and 
${\cal Y}_{red} = {\cal D}$. 
\item The family is semiuniversal at each point $(\balpha, t) \in \cM_R \times \cB_R$, that is, the Kodaira--Spencer map 
\begin{equation}
\rho:T_{\balpha, t}(\cM_R \times \cB_R) \lra H^1({\cal S}_{\balpha, t}, \Theta_{{\cal S}_{\balpha, t}}(- \log {\cal D}_{\balpha, t}))
\end{equation}
is an isomorphism. For a Zariski open subset of 
$\cM_R \times \cB_R$ on which the corresponding Okamoto--Painlev\'e pairs are of non-fibered type, 
one can choose the coordinate 
$t$ such that  (cf. Proposition \ref{prop:local})  
\begin{equation}\label{eq:timedirect}
 \rho(\frac{\partial }{\partial t}) \  \in \  H^1_{{\cal D}_{\balpha, t}}({\cal S}_{\balpha, t},  \Theta_{{\cal S}_{\balpha, t}}(- \log {\cal D}_{\balpha, t})) \hookrightarrow 
 H^1({\cal S}_{\balpha, t}, \Theta_{{\cal S}_{\balpha, t}}(- \log {\cal D}_{\balpha, t}))
\end{equation}  

\item  Let $M_R$ and $B_R$ denote the affine coordinate rings of $\cM_R$ and $\cB_R$ respectively so that $ \cM_R= \Spec M_R$ and $\cB_R = \Spec B_R $.  
(Note that $M_R$ and $B_R$ is obtained by some localization's of $ \C[\alpha_1, \cdots, \alpha_s] $ and $ \C[t]$ respectively.  ) 
 There exists a finite affine covering $\{ \tilde{U}_i \}_{i=1}^{l+k}$ of ${\cal S}$ relative to $\cM_R \times \cB_R$ 
such that there exists an isomorphism for each $i$
\begin{equation}\label{eq:covering}
  \tilde{U}_i \   \simeq \  
     \Spec (M_R \otimes B_R)  [x_i, y_i, \frac{1}{f_i(x_i, y_i, \balpha, t)}]  \subset   \Spec \C[ \balpha, t, x_i, y_i] \simeq \C^{s+3} \simeq \C^{12 - r}
\end{equation}
Here $f_i(x_i, y_i, \balpha, t)$ is a polynomial in 
 $(M_R \otimes B_R) [x_i, y_i ] $.   Moreover we may assume that $\cS -{\cal D}$ can be covered by $\{ \tilde{U_i} \}_{i=1}^l$.  Moreover for each $i$ the restriction of the 
 rational two form $\omega_{\cS}$ can be written as 
\begin{equation}
\omega_{{\cal S}|\tilde{U}_i} = \frac{dx_i \wedge dy_i}{f_i(x_i, y_i, \balpha, t)^{m_i}} 
\end{equation}

\item For each pair $i, j$ such that $\tilde{U}_i \cap 
\tilde{U}_j \not= \emptyset $, 
the coordinate transformation functions 
\begin{equation} \label{eq:transf}
x_i = f_{ij}(x_j, y_j, \balpha, t ), \quad y_i = g_{ij}
( x_j, y_j, \balpha, t)
\end{equation}
are rational functions in $ \C[x_j, y_j, \balpha, t] $.

\end{enumerate}
\end{Proposition}

Here we will give a sketch of the  proof of Proposition \ref{prop:global}.   (See  \cite{Sa-Te} for explicit 
constructions.)  
For  a generalized rational Okamoto--Painlev\'e pair $(S, Y)$, we see that  $S$  can be obtained as a blowing up of $\BP^2$ at (possibly infinitely near ) 9-points which lie on anti-canonical divisors.  Then one can parameterize these 9-points in a suitable way, and this leads to a special time parameter $t$ and other parameter 
$ \alpha_1, \cdots, \alpha_s $, hence we obtain  affine schemes $\cM_R$ and $\cB_R$,   Moreover we can construct  a semiuniversal family $\pi:{\cal S} \lra \cM_R \times \cB_R $  of rational surfaces by blowings--up of $\BP^2 \times \cM_R \times \cB_R$.  Moreover by these explicit constructions, we can obtain the affine coverings of the total space $ \cS $ as above.

\begin{Remark} {\rm We can construct a similar family of generalized Okamoto--Painlev\'e pairs of type $\tilde{A}_{r-1}$, $2 \leq r \leq 9$ (multiplicative type).  However as proved in \cite{T}, we see that 
$$
H^1_Y(S, \Theta_{S}(-\log Y)) = \{ 0  \}.
$$
This result implies that we can not obtain a differential 
equation from the generalized Okamoto--Painlev\'e pair of 
type $\tilde{A}_{r}$ as in the way above. } 
\end{Remark}

\medskip
\medskip

\section{From Global Deformations to Hamiltonian systems}
\label{sec:globaltoham}

In this section, we will explain how one can derive Hamiltonian systems from global deformation of generalized rational Okamoto-Painlev\'e pairs of additive type.  
Strictly speaking, we can obtain differential equations 
from certain special  deformations of generalized rational Okamoto--Painlev\'e pairs of additive types, but these equations  
are not always  Hamiltonian systems in the global algebraic 
coordinate systems given in Proposition \ref{prop:global}.  
In this section, we will clarify this point by means of symplectic structure on the open surfaces.  
For classical Okamoto--Painlev\'e 
pairs, it is known that  
these Hamiltonian systems are equivalent to the original Painlev\'e equations.    

\vspace{0.2cm}
Let $R=R(Y)$ be one of types of additive affine root 
systems appeared in Proposition \ref{prop:classf} and let 
\begin{equation}\label{eq:globalf2} 
\begin{array}{ccl}
 {\cal S} &  \hookleftarrow & {\cal D}      \\
\pi \downarrow \hspace{0.3cm} &    \swarrow &  \varphi  \\
 \cM_R \times \cB_R  &  &
\end{array} 
\end{equation}
be a global deformation of generalized Okamoto--Painlev\'e pairs of type $R$  as in Proposition \ref{prop:global}. 
The total space $\cS$ has a finite affine covering $\{ 
\tilde{U}_i \}_{i=1}^{l+k}$ such that 
\begin{equation}\label{eq:covering2}
  \tilde{U}_i   \simeq 
     \Spec (M_R \otimes B_R)  [x_i, y_i, \frac{1}{f_i(x_i, y_i, \balpha, t)}]  \subset   \Spec \C [ \balpha, t, x_i, y_i]\end{equation}
as in (\ref{eq:covering}).  Moreover, we may assume that  
$\cS - \cD $ can be covered by  $\{ \tilde{U}_i \}_{i=1}^l$, 
that is, 
$$
\cS - \cD = \cup_{i=1}^l \tilde{U}_{i}.
$$

Let us recall that the coordinate transformations in (\ref{eq:transf}) for $ \tilde{U}_i \cap \tilde{U}_j \not= \emptyset $
are given by the rational functions
\begin{equation} \label{eq:transf2}
x_i = f_{ij}(x_j, y_j, \balpha, t ), \quad y_i = g_{ij}
( x_j, y_j, \balpha, t)
\end{equation}
The Kodaira--Spencer class $\rho(\frac{\partial}{\partial t}) $ can be represented by the \v{C}ech 1-cocycles 
\begin{equation}\label{eq:KSC}
\rho(\frac{\partial}{\partial t}) = 
\{ \  \theta_{ij} = \frac{\partial f_{ij}}{\partial t} \frac{\partial}{\partial x_i } +  \frac{\partial g_{ij}}{\partial t} \frac{\partial}{\partial y_i } \in \Gamma( \tilde{U}_i \cap \tilde{U}_j, \Theta_{\cS/\cM_R \times \cB_R} (- \log \cD) )  \   \}  
\end{equation}
From (\ref{eq:timedirect}) of Proposition 
\ref{prop:global}, 
we may assume  that $\rho(\frac{\partial}{\partial t})$ is
 non-zero element of  the local cohomology group 
\begin{equation}
 H^1_{\cD_{\balpha, t}}(\cS_{\balpha, t}, \Theta_{\cS_{\balpha, t}}(- \log \cD_{\balpha, t}) ).
\end{equation}
Since the local cohomology group is 
the kernel of the natural restriction map (cf. Proposition \ref{prop:local}) 
\begin{equation}
\res: H^1 (\cS_{\balpha, t}, \Theta_{\cS_{\balpha, t}}
(- \log \cD_{\balpha, t})) \lra H^1( \cS_{\balpha, t} - \cD_{\balpha, t} , \Theta_{\cS_{\balpha, t}}(- \log 
\cD_{\balpha, t})), 
\end{equation}
the Kodaira--Spencer class
$ \rho(\frac{\partial}{\partial t}) $
is cohomologus to zero in $H^1( \cS_{\balpha, t} - \cD_{\balpha, t} , \Theta_{\cS_{\balpha, t}}(- \log 
\cD_{\balpha, t}))$. 

Since dimensions of these cohomology groups are constant as a function of $(\balpha, t)$,   by an argument using  the base change theorem, we see that for $ 1 \leq i \leq l$ there exist
regular vector fields 
\begin{equation}\label{eq:vectorfield}
\theta_i(x_i, y_i, \balpha, t)  = \eta_i(x_i, y_i, \balpha, t) \frac{\partial}{\partial x_i } 
+ \zeta_i(x_i, y_i, \balpha,t)  \frac{\partial}{\partial y_i } \ \in \ 
\Gamma(\tilde{U}_i, \Theta_{\tilde{U}_i}) 
\end{equation}
such that 
\begin{equation}\label{eq:cobound}
\theta_{ij}(x_i, y_i, \balpha, t)  =  \theta_{j}(x_j, y_j, \balpha, t) 
- \theta_{i}(x_i, y_i, \balpha, t).  
\end{equation}
Since we are working in the algebraic category,  
we can choose  $\eta_i(x_i, y_i, \balpha, t)$ and 
$\zeta(x_i, y_i, \balpha, t)$ as  rational functions in 
the variables $\balpha, t, x_i, y_i$.

As in (\ref{eq:lift}) of \S \ref{sec:KS}, we have the 
identity for $i, j$
\begin{equation}
\left( \frac{\partial }{\partial t}\right)_{j} = 
\left(\frac{\partial }{\partial t}\right)_{i} + \theta_{ij}(\balpha, t),   
\end{equation}
and hence just for $ 1 \leq i, j \leq l$, we have 
\begin{equation}
\left( \frac{\partial }{\partial t}\right)_{j} = 
\left(\frac{\partial }{\partial t}\right)_{i} + (\theta_j(x_j, y_j, \balpha, t)
- \theta_i (x_i, y_i, \balpha, t) ),  
\end{equation}
or 
\begin{equation}\label{eq:vf}
\left( \frac{\partial }{\partial t}\right)_{j} - \theta_j(x_j, y_j, \balpha, t)  = 
\left(\frac{\partial }{\partial t}\right)_{i} - \theta_i(x_i, y_i, \balpha, t). 
\end{equation}
This means that 
the vector fields 
\begin{equation}
  \{ \left( \frac{\partial }{\partial t}\right)_{i} - \theta_i(x_i, y_i, \balpha, t) \}_{1 \leq i \leq l}
\end{equation}
can be patched together and defines a global  vector field 
$$
\tilde{v}  \in \Gamma( \cS - \cD, \Theta_{\cS- \cD} ) 
$$
Note that this global vector field $\tilde{v}$ is a lift of 
$\frac{\partial}{\partial t}$ via $\pi: \cS - \cD \lra \cM_R \times \cB_R $. 

From the above argument, we have the following

\begin{Theorem}\label{thm:diffeq}
 Let $R= R(Y)$, $\cS, \cD, \cM_R \times \cB_R \ni (\balpha, t) $ be as above. Then we obtain the differential equation
 defined on $\cS - \cD$ whose restriction to 
 each affine chart $\tilde{U}_i $, $ 1 \leq i \leq l$,   
is given as 
\begin{equation}\label{eq:diffeq}
\renewcommand{\arraystretch}{2.2}
\left\{
\begin{array}{ccl}
\displaystyle{\frac{d x_i }{d t}} &  = & -  \eta_i (x_i, y_i, \balpha, t) \\
\displaystyle{\frac{d y_i }{d t}} & = & -\zeta_i (x_i, y_i. \balpha, t) 
\end{array}
\right.
\end{equation}
where the functions appeared in the  right hand sides are 
rational functions in the variables $x_i, y_i, \balpha, t$. 

\end{Theorem}

\begin{Remark}{\rm 
\begin{enumerate}

\item The argument above shows that 
there exists a differential equation as above at least locally for any direction
corresponding to the kernel of the restriction map
$$
\res: H^1(S, \Theta_S( -\log D)) \lra H^1(S -D, \Theta_S( - \log D)). 
$$
\item 
Let us recall  the so-called Painlev\'e property which  is states as follows. {\em If $(x(t), y(t))$ is a local solution of (\ref{eq:diffeq}) determined by an arbitrary initial conditions $(x_0 = x(t_0), y_0=y(t_0)) \in \tilde{U}_i $ with fixed $t_0 \in \cB_R $ then both solutions $x(t)$ and $y(t)$ can be meromorphically continued along any curve in $\cB_R$ with a starting point $t_0$.} For non-classical Okamoto--Painlev\'e pair, it is not clear that the differential equation 
in (\ref{eq:diffeq}) has the Painlev\'e property.  In general, the proof of Painlev\'e property  for classical 
Painlev\'e equation is not so easy.  We hope that there is an easy geometric proof of the Painlev\'e property for differential equation in (\ref{eq:diffeq}).
\end{enumerate}
}
\end{Remark}

It is well-known that each  classical Painlev\'e differential equation $P_J$, $J =I, II, \cdots, VI$ 
 is equivalent to a Hamiltonian system $(H_J)$ whose 
 Hamiltonian $ H_J(x, y, \balpha, t) $ is 
a polynomial in $ (x, y) \in \C^2 $ (cf. \cite{O1}, \cite{MMT}).  

\begin{equation}
\renewcommand{\arraystretch}{2.2}
(H_J): \quad
\left\{
\begin{array}{ccl}
\displaystyle{ \frac{d x}{d t}} & = & 
\displaystyle{ \frac{ \partial H_J}{ \partial y } } \\
\displaystyle{ \frac{ d y}{d t}} & = & - 
\displaystyle{ \frac{ \partial H_J}{\partial x} }.
\end{array}
\right. 
\end{equation}
\vspace{0.2cm}

In what follows, we shall show how this Hamiltonian systems arise from our differential equations in (\ref{eq:diffeq})  obtained from the deformation of generalized rational 
Okamoto--Painlv\'e pairs.   

Let us recall the general situation. Recall  that 
 $\cS - \cD$ is covered by  $\{ \tilde{U}_i \}_{i=1}^l $ and 
$$
\tilde{U}_{i} =    \Spec (M_R \otimes B_R)  [x_i, y_i, \frac{1}{f_i( x_i, y_i, \balpha, t) } ] \subset \cM_R \times \cB_R \times \C^2, 
$$
and the restriction of the relative two form $\omega_{\cS}$ to $\tilde{U}_i$  can be written as 
\begin{equation}\label{eq:form}
\omega_{ {\cal S} |\tilde{U}_i } = \frac{dx_i \wedge d y_i}{f_i( x_i, y_i, \balpha, t)^m}
\end{equation}
Let $ \theta_i(x_i, y_i, \balpha, t) $ be the vector fields 
defined in (\ref{eq:vectorfield}).  
The contraction of $\theta_i$ and $\omega_{{\cal S} |\tilde{U}_i}$ is given by 
$$
 \theta_i \cdot \omega_{\cS|\tilde{U}_i} = \frac{1}{f_i^{m_i}}( \eta_i d y_i - \zeta_i d x_i ).   
$$
Consider the following  regular two form on $\tilde{U}_i$ for each $ 1 \leq i \leq l $
$$
\Omega_{i}:=  \omega_{\cS|\tilde{U}_i} -  
( \theta_i \cdot \omega_{\cS|\tilde{U}_i} )  \wedge dt. 
$$

\begin{Lemma}\label{lem:ham}
On $\tilde{U}_i \cap \tilde{U}_j \not= 
\emptyset $, we have 
$$
\Omega_i = \Omega_j \in \Gamma(\tilde{U}_i \cap \tilde{U}_j, 
\Omega^2_{\cS-\cD/\cM_R}).  
$$
Hence, we have a regular two form 
$ \Omega \in \Gamma(\cS - \cD, \Omega^2_{\cS- \cD/\cM_R})$ such that
$$
\Omega_{|\tilde{U}_i} = \Omega_i
$$
\end{Lemma}

{\it Proof.}
Since $\pi:\cS - \cD \lra \cM_R \times \cB_R$ is smooth,  we have the following exact sequence 
$$
0 \lra \pi^*\Omega^1_{\cM_R \times \cB_R/\cM_R} \lra 
\Omega^1_{\cS- \cD /\cM_R} \lra \Omega^1_{\cS- \cD/\cM_R \times \cB_R} \lra 0.
$$
Moreover since  the 
relative dimension of $\pi:\cS- \cD \lra \cM_R \times \cB_R$ 
is two, we have the exact sequence 
$$
0 \lra \Omega^1_{\cS - \cD/\cM_R \times \cB_R} \otimes \pi^*\Omega^1_{\cM_R \times \cB_R/\cM_R} \lra \Omega^2_{\cS- \cD/\cM_R} \lra 
\Omega^2_{\cS- \cD/\cM_R \times \cB_R} \lra 0
$$
Note that  the global section $\omega_{\cS}$ lies in the space  
$$
\omega_{\cS} \in \Gamma(\cS - \cD, \Omega^2_{\cS- \cD/\cM_R \times \cB_R}) 
$$
Hence by a local calculation  if we restrict $\omega_{\cS}$ to each $\tilde{U}_{i}$, then  on $\tilde{U}_i \cap \tilde{U}_j \not= \emptyset $ we have the relation 
$$
\omega_{\cS|\tilde{U}_i} = \omega_{\cS|\tilde{U}_j} 
- \theta_{ij} \cdot \omega_{\cS|\tilde{U}_i} \wedge dt
$$
where $\theta_{ij}$ is Kodaira--Spencer class representing 
$\rho(\frac{\partial}{\partial t})$.  Then,  
by using the relation (\ref{eq:cobound}), 
we see that 
\begin{equation}
\omega_{\cS|\tilde{U}_i}- \theta_i \cdot \omega_{\cS|\tilde{U}_i} \wedge dt = \omega_{\cS|\tilde{U}_j} 
- \theta_{j} \cdot \omega_{\cS|\tilde{U}_j} \wedge dt. 
\end{equation}
This completes the proof. \qed

Let 
$$
d_{\cS-\cD/\cM_R}:
\Omega^2_{\cS -\cD/ \cM_R} \lra  \Omega^3_{\cS -\cD/ \cM_R}
$$
be the relative exterior derivative. 
Since the deformation $\cS- \cD \lra \cB_R$ preserves the 
regular two form $\omega_{\cS_{\balpha, t}}$, by an standard argument we have the following
\begin{Proposition} \label{prop:closed}
$$
d_{\cS-\cD/\cM_R}(\Omega) = 0
$$
\end{Proposition}
Looking at the isomorphism 
$$
\Omega^3_{\cS -\cD/\cM_R} \simeq \Omega^2_{\cS - \cD/\cM_R \times \cB_R} \otimes \pi^*(\Omega^1_{\cM_R \times \cB_R/\cM_R}),
$$
let us write
$$
d_{\cS-\cD/\cM_R}(\Omega) = \eta_{\cS - \cD} \wedge dt
$$
where
$$
\eta_{\cS - \cD} \in \Gamma(\cS -\cD, \Omega^2_{\cS - \cD/\cM_R \times \cB_R}).  
$$
(Note that $\eta_{\cS -\cD}$ may not be  a global regular 
2-form in $\Omega_{\cS - \cD/\cM_R}$. )
 Then we have 
 \begin{eqnarray*}
 d_{\cS - \cD/\cM_R} (\Omega_{|\tilde{U}_i}) & = &  
 d_{\cS - \cD/\cM_R}(\frac{dx_i \wedge dy_i}{f_i(x_i, y_i, \balpha,t)^{m_i}} - 
 (\theta_i \cdot  \omega_{\cS}) \wedge dt) \\  
  & = &   ( \frac{\partial}{\partial t}\left( \frac{1}{f_i(x_i, y_i, \balpha, t)^{m_i}} \right ) dx_i \wedge dy_i - d_{\pi} (\theta_i \cdot \omega_{\cS} ))\wedge dt 
  \end{eqnarray*}
 where we set $d_{\pi} = d_{\cS - \cD/\cM_R \times \cB_R}$. 

Therefore, 
Proposition \ref{prop:closed} implies the following important

\begin{Corollary}
  For each $i$, $1 \leq i \leq l$,  we have the fundamental equation 

\begin{equation}\label{eq:fund}
\frac{\partial}{\partial t}\left( \frac{1}{f_i(x_i, y_i, \balpha, t)^{m_i}} \right ) dx_i \wedge dy_i - d_{\pi} (\theta_i \cdot \omega_{\cS} ) = 0 
\end{equation}

\end{Corollary}

\vspace{0.5cm}

Now we obtain the following fundamental results. 

\begin{Proposition}\label{prop:hamiltonian}
For $i$, $1 \leq i \leq l$ such that 
\begin{equation}
\tilde{U}_i =
 \cM_R \times \cB_R \times \Spec \C[x_i, y_i]\simeq 
 \cM_R \times \cB_R \times \C^2
, \quad \omega_{\cS|\tilde{U}_i} = dx_i \wedge dy_i, 
\end{equation}
we have 
$$
d_{\pi}(\theta_i \cdot dx_i \wedge dy_i) = d_{\pi}( \eta_i d y_i - \zeta_i d x_i)  = 0.
$$
Since $H^1_{DR}(\C^2) = 0$, we have a regular function 
$H_i(x_i, y_i, \balpha, t) \in (M_R \otimes B_R)[x_i, y_i]$ such that 
$$
d_{\pi} H_i =  \frac{ \partial H_i}{\partial x_i} dx_i + 
 \frac{ \partial H_i}{\partial y_i} dy_i =  - (\theta_i \cdot dx_i \wedge dy_i) =- \eta_i d y_i + \zeta_i d x_i .
$$
From this, we have 
$$
- \eta_i =   \frac{ \partial H_i}{\partial y_i}, - \zeta_i = -
 \frac{ \partial H_i}{\partial x_i}.  
$$
Therefore, the differential equation (\ref{eq:diffeq}) can be written in the  Hamiltonian system
\begin{equation}
\renewcommand{\arraystretch}{2.2}
\left\{
\begin{array}{ccl}
\displaystyle{ \frac{d x_i}{d t}} & = & 
\displaystyle{ \frac{ \partial H_i}{ \partial y_i } } \\
\displaystyle{ \frac{ d y_i}{d t}} & = & - 
\displaystyle{ \frac{ \partial H_i}{\partial x_i} }.
\end{array}
\right. 
\end{equation}
\end{Proposition}

\begin{Remark}
{\rm \begin{enumerate}
\item If $f_i(x_i, y_i, \balpha, t)$ in (\ref{eq:form})  is independent of $t$, from the fundamental equation (\ref{eq:fund}),we obtain
$$
d_{\pi}( \theta_i \cdot \frac{dx_i \wedge dy_i}{f_i(x_i, y_i, \balpha)^{m_i}}) = 0.
$$  
Therefore, we may have a chance to have a regular function 
$ H_i(x_i, y_i, \balpha, t) $ on $\tilde{U}_i$ such that 
$$
d_{\pi} H_i =  - \theta_i \cdot \frac{dx_i \wedge dy_i}{f_i(x_i, y_i, \balpha)^{m_i}}.
$$
In this case, the differential equation in (\ref{eq:diffeq})  can be written in 
\begin{equation}
\renewcommand{\arraystretch}{2.2}
\left\{
\begin{array}{ccl}
\displaystyle{ \frac{d x_i}{d t}} & = & 
\displaystyle{ (f_i)^{m_i} \cdot \frac{ \partial H_i}{ \partial y_i } } \\
\displaystyle{ \frac{ d y_i}{d t}} & = & - 
\displaystyle{ (f_i)^{m_i} \cdot \frac{ \partial H_i}{\partial x_i} }, 
\end{array}
\right. 
\end{equation}
or equivalently, 
\begin{equation}
\renewcommand{\arraystretch}{2.2}
\left\{
\begin{array}{ccl}
\displaystyle{\frac{1}{(f_i)^{m_i}} \frac{d x_i}{d t}} & = & 
\displaystyle{   \frac{ \partial H_i}{ \partial y_i } } \\
\displaystyle{ \frac{1}{(f_i)^{m_i}} \frac{ d y_i}{d t}} & = & - 
\displaystyle{  \frac{ \partial H_i}{\partial x_i} }.
\end{array}
\right. 
\end{equation}
\item In general, we can not transform the differential equation in (\ref{eq:diffeq}) into a Hamiltonian system in the global 
affine coordinates.  

\item Takano, et al. show that for any  
Okamoto--Painlev\'e pair $(S, Y)$ 
of type $\tilde{D}_4 (=P_{VI}), 
\tilde{D}_5 (=P_{V}), \tilde{D}_6 (=P_{III}), \tilde{E}_6 (=P_{IV}), \tilde{E}_7 (= P_{II})$ , the open surface 
$S - Y_{red}$ is covered by a finite number of affine spaces  $ {U}_i = \C^2$ and regular 2-form $\omega_{\cS|\tilde{U}_i} $ can be written as in $ dx_i \wedge dy_i$.  
Hence from Proposition \ref{prop:hamiltonian}  we obtain the Hamiltonian systems for those Okamoto--Painlev\'e pairs 
on any affine chart $U_i$ of $ S - D $ as proved in \cite{O1}, \cite{MMT}.  Note that for an Okamoto--Painlev\'e pair $(S, Y)$ of type $ \tilde{D}_8$, $S - Y_{red}$ does not contain $\C^2$  (cf. Theorem \ref{thm:classf} and 
Proposition \ref{prop:classf}). For 
explicit descriptions of  $\tilde{E}_7$ and 
$\tilde{D}_8$, see \S \ref{sec:examples}.  
\end{enumerate}
}
\end{Remark}

We summarize our results in this section as follows (cf. Figure \ref{fig:fund}).  

\begin{figure}[h]
\begin{tabular}{c}
\fbox{
\begin{tabular}{l} 
Deformation $  \cD \hookrightarrow  \cS \lra \cM_R \times \cB_R \ni (\balpha, t)  $  of Okamoto--Painlev\'e pairs   
such that \\   for any $ (\balpha, t)  \in \cB_R $
the Kodaira--Spencer class $\rho(\frac{\partial}{\partial t})$ lies in the Kernel of the \\  restriction map  
$
  \res: H^1 (\cS_{\balpha, t}, \Theta_{\cS_{\balpha, t}}
(- \log \cD_{\balpha, t})) \lra H^1( \cS_{\balpha, t} - \cD_{\balpha, t} , \Theta_{\cS_{\balpha, t}}(- \log 
\cD_{\balpha, t})) 
$
\end{tabular} }  \\
 \\
$\Downarrow$ \\
   \\
\fbox{There exists a global holomorphic vector field $ \tilde{v}$ on $\cS - \cD $ which is a lift of $\frac{\partial}{\partial t}$ } \\
  \\
$\Downarrow$ \\
  \\
\fbox{
Differential Equations on $\pi:\cS - \cD \lra \cM_R \times \cB_R  $ defined by $\tilde{v}$ } \\
  \\
\quad \hspace{3cm} $\Downarrow$ \quad \fbox{Painlev\'e property} \\
  \\
\fbox{Local trivializations of the deformation $\cS - \cD  \lra  \cB_R $ }.
\end{tabular}
\vspace{0.3cm}
\caption{}
\label{fig:fund}
\end{figure}

\newpage
\vspace{1cm}
\section{Painlev\'e Equations}\label{sec:painleve}

Let us recall the classical Painlev\'e differential equations and Hamiltonian systems which are equivalent to the Painlev\'e equations (\cite{IKSY}, \cite{T}, \cite{O1}). 

Painlev\'e equations $P_J$ $(J = I, II, \cdots, VI)$ are given in Table \ref{tab:painleve}:

\begin{table}[h]
\begin{center}
\vspace{0.5cm}
$$
\renewcommand{\arraystretch}{2.6}
\begin{array}{cccl}
P_{I}: & \displaystyle{\frac{d^2 x}{d t^2}} &= & 6 x^2 + t, \\
P_{II}: & \displaystyle{\frac{d^2 x}{d t^2}} & =&  2 x^3 + t x +  \alpha,  \\
P_{III}: & \displaystyle{\frac{d^2 x}{d t^2}}  & = & 
\displaystyle{
\frac{1}{x} \left( \frac{d x}{d t} \right)^2 - \frac{1}{t} \frac{d x}{d t} + \frac{1}{t} (\alpha x^2 + \beta) + \gamma x^3 + \frac{\delta}{x},  }  \\
P_{IV}: & \displaystyle{\frac{d^2 x}{d t^2}}  & = & \displaystyle{
\frac{1}{2x} \left( \frac{d x}{d t} \right)^2 + \frac{3}{2} x^3 + 4 t x^2 + 2(t^2 - \alpha) x + \frac{\beta}{x}},      \\
P_{V}: & \displaystyle{\frac{d^2 x}{d t^2}} & = & 
\displaystyle{\left(\frac{1}{2x} + \frac{1}{x-1} \right) \left( \frac{d x}{d t} \right)^2 -  \frac{1}{t} \frac{d x}{d t} + \frac{(x-1)^2}{t^2} \left( \alpha x + \frac{\beta}{x} \right) + \gamma \frac{x}{t} + \delta \frac{x(x+1)}{x-1}   },  \\
P_{VI}: & \displaystyle{\frac{d^2 x}{d t^2}}  &= & 
\displaystyle{ \frac{1}{2}\left( \frac{1}{x} + \frac{1}{x-1}+ \frac{1}{x-t} \right)\left( \frac{d x}{d t} \right)^2 -  \left( \frac{1}{t} + \frac{1}{t-1}+ \frac{1}{x-t} \right)\left( \frac{d x}{d t} \right)},     \\
        &  & &   \displaystyle{ + \frac{x(x-1)(x-t)}{t^2(t-1)^2}
         \left[ \alpha - \beta \frac{t}{x^2} + \gamma 
         \frac{t-1}{(x-1)^2} + \left(\frac{1}{2} - \delta \right) \frac{t(t-1)}{(x-t)^2} \right] }. 
\end{array}
$$
\end{center}
\vspace{0.5cm}
\caption{}
\label{tab:painleve}
\end{table}

Here $x$ and $t$ are complex variables, 
$ \alpha, \beta, \gamma $ and $\delta$ are complex constants. It is known that each $P_J$ is equivalent to a Hamiltonian 
system (cf. \cite{O1}, \cite{IKSY}, \cite{MMT}): 

\begin{equation}\label{eq:hamiltonian}
\renewcommand{\arraystretch}{2.2}
(H_J): \left\{
\begin{array}{ccl}
\displaystyle{ \frac{d x}{d t}} & = & 
\displaystyle{ \frac{ \partial H_J}{ \partial y} },  \\
\displaystyle{ \frac{ d y}{d t}} & = & - 
\displaystyle{ \frac{ \partial H_J}{\partial x} }, 
\end{array}
\right. 
\end{equation}
where the Hamiltonians $H_J$ are given in Table \ref{tab:ham}.   

\begin{table}[h]
$$
\renewcommand{\arraystretch}{2.2}
\begin{array}{ccl}
H_{I}(x, y, t) & = &\displaystyle{ \frac{1}{2} y^2 -2 x^3 - t x }, \\
H_{II}(x, y, t) & = & \displaystyle{\frac{1}{2} y^2 - \left( x^2 + \frac{t}{2} \right) y - \left( \alpha + \frac{1}{2} \right) x}, \\
H_{III}(x, y, t) & = & \displaystyle{\frac{1}{t} \left[
2 x^2 y^2 - \left\{2 \eta_{\infty} t x^2 + (2 \kappa_0 +1) x - 2 \eta_0 t  \right\} y + \eta_{\infty} \left( \kappa_0 + \kappa_{\infty}\right ) t x \right]}, \\
H_{IV}(x, y, t) & = & \displaystyle{2 x y^2 - \left\{ x^2 + 
2 t x + 2 \kappa_0  \right \} y + \kappa_0 x}, \\
H_{V}(x, y, t) & = & \displaystyle{
\frac{1}{t} \left[x(x-1)^2 y^2 - \left\{\kappa_0 (x-1)^2 + 
\kappa_t x(x-1) - \eta t x  \right\} y +  \kappa (x-1) \right] }, \\
  & & \displaystyle{ \left( \kappa := \frac{1}{4} \left\{ (\kappa_0 + \kappa_t)^2 - \kappa^2_{\infty} \right\} \right)}, \\
  H_{VI}(x, y, t) & = & \displaystyle{ \frac{1}{t(t-1)} 
  \left[ x (x-1)(x-t) y^2 - 
  \left\{ \kappa_0 (x -1) (x - t) \right. \right.} \\  
      & & \displaystyle{ \left. \left.  + \kappa_1 x(x-t) + (\kappa_t - 1) x (x -1) 
      \right\} y + \kappa (x - t) \right]} \\
      & &  \displaystyle{ \left( \kappa := \frac{1}{4} 
      \left\{ (\kappa_0 + \kappa_1 +  \kappa_t -1)^2 - \kappa^2_{\infty} \right\} \right) }. 
\end{array}
$$
\caption{}
\label{tab:ham}
\end{table}

Moreover the relations between the constants in the equations 
$P_J$ and the Hamiltonians $H_J$ are given in Table \ref{tab:auxparm}. 

\begin{table}[h]
\begin{center}
\renewcommand{\arraystretch}{1.7}
\begin{tabular}{|c||c|c|c|c|c|} \hline
        & $\alpha$ & $\beta$ & $\gamma$ & $\delta$ & number of aux. parameters \\  \hline 
      $  P_{I}$ & none & none   & none   & none   & 0 \\ \hline 
       $ P_{II}$ & $ \alpha $ & none  & none  & none & 1 \\ \hline 
$P_{III}$ & $ - 4 \eta_{\infty} \kappa_{\infty} $ & 
$ 4 \eta_{\infty}(\kappa_0 +1)$ & $4 \eta^2_{\infty}$ & $-4 \eta^2_{0}$ & 2 \\ \hline
$P_{IV}$  &  $- \kappa_0 + 2 \kappa_{\infty} + 1$ & $-2 \kappa^2_0$ & none  & none  & 2  \\ \hline 
$P_{V}$ & $\frac{1}{2} \kappa^2_{\infty}$ & $ -\frac{1}{2} 
\kappa^2_0$ & $ - \eta(1 + \kappa_t ) $& $-\frac{1}{2}\eta^2/2 $  & 3 \\ \hline 
$P_{VI}$ &   $\frac{1}{2} \kappa^2_{\infty}$ & $ \frac{1}{2} 
\kappa^2_0 $ &  $ \frac{1}{2} \kappa_1^2 $  & $ \frac{1}{2} \kappa^2_t $ & 4 
\\ \hline 
\end{tabular}
\end{center}
\vspace{0.5cm}
\caption{}
\label{tab:auxparm}
\end{table}

For the meaning of the constants in Table \ref{tab:auxparm}, see \cite{IKSY}, \cite{O2}.  Note that these constants are not effective parameters.   In some cases, we can normalize these constants further by coordinate transformations.  
Moreover, the equivalence of $P_J$ and $(H_J)$ means that 
if we eliminate the variable $y$ in $(H_J)$ then we obtain $(P_J)$.

\newpage
\begin{Remark}{\rm 
For the Painlev\'e equation of type $P_{III}$, we have the following remark by Sakai in \cite{Sakai}.  
The Painlev\'e equation of type $P_{III}$ as in Table \ref{tab:painleve} can be 
transformed into 
\begin{equation}\label{eq:PIII}
\frac{d^2 x}{d t^2}  =  \frac{1}{x} \left( \frac{d x}{d t} \right)^2 - \frac{1}{t} \frac{d x}{d t} + \frac{x^2}{4t^2} (\gamma x + \alpha) +\frac{\beta}{4t}+  \frac{\delta}{4x}. 
\end{equation}
If $ \gamma \delta \not= 0 $, then we can normalize $ \gamma = - \delta = 4 $ without loss of generality.  In this case we obtain the Painlev\'e equation of type
 $ P^{\tilde{D_6}}_{III} $:

\begin{equation}
P^{\tilde{D_6}}_{III}: \quad \frac{d^2 x}{d t^2}  =  \frac{1}{x} \left( \frac{d x}{d t} \right)^2 - \frac{1}{t} \frac{d x}{d t} + \frac{x^2}{4t^2} ( 4x + \alpha) + \frac{\beta}{4t}-  \frac{1}{x}. 
\end{equation}

If one of $\gamma $ and $\delta$ equals to zero (not both), then we have $P^{\tilde{D_7}}_{III}$.  
\begin{equation}\label{eq:difd7}
P^{\tilde{D_7}}_{III}:  \quad \frac{d^2 x}{d t^2}  =  \frac{1}{x} \left( \frac{d x}{d t} \right)^2 - \frac{1}{t} \frac{d x}{d t} - \frac{4 x^2}{t^2} - \frac{1+2a}{t}.
\end{equation}

Moreover 
when $\gamma = \delta = 0$, we have $P^{\tilde{D_8}}_{III}$:
\begin{equation}\label{eq:difd8}
P^{\tilde{D_8}}_{III}: \quad \frac{d^2 x}{d t^2}  =  \frac{1}{x} \left( \frac{d x}{d t} \right)^2 - \frac{1}{t} \frac{d x}{d t} - \frac{x^2}{4t^2} - \frac{4}{t}.
\end{equation}
These differential equations correspond to 
generalized rational Okamoto--Pailev\'e pairs of type 
$\tilde{D_6}, \tilde{D_7}, \tilde{D_8}$ respectively.  
}
\end{Remark}

\newpage
\section{Examples}
\label{sec:examples}

In this section, we will apply our methods for deriving 
the differential equation in (\ref{eq:diffeq}) from the 
explicit deformations of Okamoto--Painlev\'e pairs.  
We shall give a full detail of  the cases of  $\tilde{E_7}(= P_{II})$  and 
$ \tilde{D_8}(= P^{\tilde{D_8}}_{III}) $.  
For other cases, see \cite{Sa-Te}.  

\vspace{0.2cm}
\begin{Example}\label{eq:E7}{\rm 
{\bf $\tilde{\bold E}_7$--type:}  In this case, we will use the Takano's 
description of  the family of   Okaomoto--Painlev\'e pairs of type $\tilde{E_7}$ (cf. [Theorem 4, \cite{MMT}]). 
We will not consider  all of the family $\cS  \lra \cM_R \times \cB_R$, but  consider the family $\cS - \cD \lra \cM_R \times \cB_R$ which is constructed as follows  (cf. \cite{SU}).  
Let us set 
$$
\cM_R = \Spec \C[\alpha] \simeq \C, \quad \cB_R = \Spec \C[t] \simeq \C   
$$
and take three affine schemes $i =1, 2, 3$
\begin{equation}
\tilde{U}_{i} = \Spec \C[\alpha, t, x_i, y_i] \simeq \C^4.
\end{equation}
The family $\cS - \cD \lra \cM_R \times \cB_R$ can be constructed  by
patching these affine schemes by 
the coordinate transformations given as follows (cf. \cite{MMT}):
$$
\renewcommand{\arraystretch}{1.3}
\begin{array}{lll}
x_0 & = \displaystyle{ \frac{1}{x_1}} & = \displaystyle{\frac{1}{x_2}} \\
y_0 & = x_1 ((- \alpha-\frac{1}{2}) - x_1y_1) & =2 x_2^{-2} + t + (\alpha-\frac{1}{2})x_2 - y_2 x_2^2  \\
x_1 & = x_2 & = \displaystyle{\frac{1}{x_0}} \\  
 y_1 &  = \displaystyle{- \frac{2}{x_2^4} - \frac{t}{x_2^2} - \frac{2 \alpha}{x_2} + y_2} & = \displaystyle{x_0((-\alpha- \frac{1}{2}) - x_0 y_0) }  \\
x_2 & = \displaystyle{\frac{1}{x_0}} & = x_1 \\
y_2 & = \displaystyle{2 x_0^4 + t x_0^2 + (\alpha-\frac{1}{2}) x_0 - x_0^2 y_0} & = \displaystyle{\frac{2}{x_1^{4}} + \frac{t}{x_1^2} + \frac{2 \alpha}{x_1} + y_1 } 
\end{array}.
$$
The Kodaira--Spencer class $\rho(\frac{\partial}{\partial t})_{\cS_t - \cD_t} $ is given by the \v{C}ech 1-cocycle
\begin{equation}
\theta_{01} = 0, \quad \theta_{02} = \frac{\partial}{\partial y_0}, \quad \theta_{12} = -x_1^{-2} \frac{\partial}{\partial y_1}.
\end{equation}
Setting 
\begin{eqnarray}
\theta_0 & := & \left[- y_0 +x_0^2 + \frac{t}{2} \right] \frac{\partial}{\partial x_0}- \left[ 2x_0 y_0+ \alpha + \frac{1}{2} \right] \frac{\partial}{\partial y_0} \\ 
\theta_1 & := & \frac{1}{2}\left[-2 - tx_1^2 - x_1^3 - 2 \alpha x_1^3 - 2 x_1^4 y_1 \right]  \frac{\partial}{\partial x_1}  \\
   & & \quad  \quad + \frac{1}{4}\left[(1 + 2 \alpha + 4 x_1 y_1 )(t + x_1 (1 + 2 \alpha + 2 x_1 y_1))\right] \frac{\partial}{\partial y_1}  \\
\theta_2 & := &  \frac{1}{2} ( 2 + t x_2^2 + (2 \alpha -1) x_2^3 - 2 x_2^4 y_2) \frac{\partial}{\partial x_2} \nonumber \\
  &  & \quad \quad  + \frac{1}{4} ( -1 + 2 \alpha-4 x_2 y_2)(t + x_2(-1 + 2\alpha - 2 x_2 y_2)) \frac{\partial}{\partial y_2},
\end{eqnarray}
we have the relations
$$
\theta_{01} = \theta_{1}  - \theta_{0}, \quad  \theta_{02} = \theta_2 - \theta_0, \quad \theta_{12} = \theta_2 - \theta_1,
$$
that is, 
$$
\theta_0 = \theta_1, \quad \theta_2 = \theta_{02} + \theta_0.
$$
Since  on each $\tilde{U_i}$, the relative 2-forms $\omega_{\cS - \cD}$ is given by 
$$
\omega_{\cS - \cD| \tilde{U}_i} = d x_i \wedge d y_i, 
$$
by Proposition  \ref{prop:hamiltonian} the 1-forms $ \theta_i dx_i \wedge dy_i $ is 
exact form, hence there exists a polynomial $ H_i(x_i, y_i, \alpha, t) $ satisfying 
$$
 - \theta_i dx_i \wedge dy_i = d_{\pi} H_i.  
$$
The polynomials $H_i $ are called the Hamiltonians and given by
\begin{eqnarray}
\renewcommand{\arraystretch}{1.3}
\label{eq:ham2}
H_0(x_0, y_0, \alpha, t) & = & \displaystyle{\frac{1}{2} y_0^2 - \left( x_0^2 + \frac{t}{2} \right) y_0 - \left( \alpha + \frac{1}{2}\right) x_0},  \\
H_1(x_1, y_1, \alpha, t) & = & \displaystyle{\frac{t x_1}{4} + \frac{
\alpha t x_1}{2} + \frac{x_1^2}{8} + \frac{\alpha x_1^2}{2} + \frac{\alpha^2 x_1^2}{2}+ y_1 + \frac{1}{2} t x_1^2 y_1} \\
& & \quad \quad  \displaystyle{ + 
\frac{x_1^3 y_1}{2} + 
    \alpha  x_1^3 y_1 + \frac{x_1^4 y_1^2}{2} },  \nonumber \\
H_2(x_2, y_2, \alpha, t) & = & \displaystyle{\frac{1}{8}\left[(1 - 2 \alpha)^2 x_2^2 - 8 y_2 -  
        4(-1 + 2\alpha) x_2^3 y_2 + 4 x_2^4 y_2^2 \right.} \\
       &  &  \quad \quad \displaystyle{ \left.  -  2 t x_2(1 - 2 \alpha + 2x_2 y_2) \right)]}.
\end{eqnarray}
Hence the Hamiltonian system defined on $\tilde{U}_0$ is given by 
\begin{equation}
\label{eq:hamsys2}
\renewcommand{\arraystretch}{2.0}
\left\{
\begin{array}{ccl}
 \displaystyle{\frac{d x_0}{dt}} & = &
\displaystyle{ \frac{\partial H_0}{\partial y_0 } = y_0 - x_0^2 - \frac{t}{2} },  \\
\displaystyle{ \frac{d y_0}{d t}} & = & \displaystyle{ - \frac{\partial H_0}{\partial x_0 } = 2 x_0 y_0+ \alpha + \frac{1}{2} }.
\end{array}  \right. 
\end{equation}
Eliminating $y_0$ in (\ref{eq:hamsys2}), we obtain 
\begin{equation}
\frac{d^2 x_0}{d t^2} = 2 x_0^3 + x_0 t +  \alpha, 
\end{equation}
which is  the Painlev\'e equation $P_{II}$ in Table \ref{tab:painleve}.  
}
\end{Example}

\vspace{0.2cm}

\begin{Example} {\bf $\tilde{\bold D}_8$--type}: {\rm The Okamoto--Painlev\'e pair $(S, Y)$ 
of type $ \tilde{D}_8 $ did not appear in the former literatures  \cite{IKSY}, \cite{O1} explicitly.  Since $S -Y_{red}$ does not contain $ \C^2 $ as a Zariski open set ( cf. Theorem \ref{thm:classf}), the situation 
is a little bit different from the classical cases. 

We can construct a family of generalized rational 
Okamoto--Painlev\'e pair of type $ \tilde{D}_8 $  $ \pi: \cS - \cD \lra \cB_R $ by blowings-ups as in Sakai \cite{Sakai}.  For detail, see \cite{Sa-Te}.  
Here  note that $ \dim H^1(\cS_t, \Theta_{\cS_t}( - \log \cD_t)) = 1 $ and $ \dim \cM_R = 0 $ and  
$$
\cB_R = \Spec \C[ t, t^{-1}] \simeq \C^{\times}. 
$$

The total space $ \cS - \cD $ is covered by the three 
affine open sets:
$$ 
\cS- \cD \ = \tilde{U}_0 \cup \tilde{U}_1 \cup \tilde{U}_2.
$$
These affine open sets are given by: 
\begin{eqnarray}
	\tilde{U}_0 & = & \Spec \C [x_0,y_0,\frac{1}{y_0}, t, t^{-1} ] \ \cong
  ({\C}^2-\{y_0=0\}) \times \C^{\times},  \nonumber\\
\tilde{U}_1& = & \Spec \C [x_1,y_1,\frac{1}{F_1(x_1,y_1)}, t, t^{-1} ] \ \cong \ ({\C}^2-\{ F_1(x_1,y_1) =0\}) \times 
\C^{\times},  \nonumber\\
	\tilde{U}_2 & = & 
	\Spec \C [x_2,y_2,\frac{1}{F_2(x_2,y_2,t)}, t, t^{-1}] \ \cong \  {\C}^3-\{F_2(x_2,y_2,t)=0, \ t = 0\}, \nonumber
\end{eqnarray}
where
\begin{eqnarray}
	F_1(x_1,y_1) & = & 1+x_1 {y_1}^2,  \\
	F_2(x_2,y_2,t) & = & t-t y_2+x_2 {y_2}^2. 
\end{eqnarray}

The coordinate transformations are given  as follows
$$
\renewcommand{\arraystretch}{2.3}
\begin{array}{cll}
x_0 & =  \displaystyle{\frac{1}{y_1F_1(x_1,y_1)}} &  =\displaystyle{\frac{1}{y_2}},  \\ 
y_0 & =\displaystyle{\frac{1}{{y_1}^2 F_1(x_1,y_1)}}& = {y_2}^2\ (t+x_2  {y_2}^2-ty_2 ),   \\
x_1 & =\displaystyle{\frac{{y_0}^2  (-{x_0}^2+y_0)}{{x_0}^4}}  & = y_2^6 (t + x_2 y_2^2 - t y_2),   \\
y_1 & = \displaystyle{\frac{x_0}{y_0}} & = \displaystyle{\frac{1}{y_2^3 ( t( 1 - y_2) + x_2 y_2^2 )} },  \\
x_2 & = x_0  (-t x_0+ x_0^3  y_0+t) & = 
\displaystyle{\frac{1+ t y_1^4 F_1^3 (y_1^2 F_1 - 1)}{y_1^6 F_1^5}},   \\
y_2 & =\displaystyle{\frac{1}{x_0}} &  = y_1^2(1 + x_1 y_1^2).  
\end{array}
$$

The Kodaira--Spencer class 
$ \rho(\frac{\partial}{\partial t})_{\cS_t - \cD_t} \in  H^1(\cS_t  - \cD_t ,\Theta_{\cS_t}(-\log \cD_t)) $
is given by \v{C}ech cocycles
\begin{equation}
{\theta }_{01} = 0, \quad \displaystyle{{\theta }_{02} = \frac{-1+x_0}{{x_0}^3}\frac{\partial }{\partial y_0}}, 
\quad 	\displaystyle{{\theta }_{21} = \frac{-1+y_2}{{y_2}^2}\frac{\partial }{\partial x_2}}. 
\end{equation}

Since $\rho(\frac{\partial}{\partial t})_{\cS_t - \cD_t}  = 0 \in  H^1(\cS_t  - \cD_t ,\Theta_{\cS_t}(-\log \cD_t)) $, we 
can obtain \v{C}ech coboundary $\{ \theta_i \in \Gamma( \tilde{U}_i, \Theta_{\tilde{U}_i/\cB_R }  
\} $ such that 
$$ 
\{ \theta_{ij}\} =\{ \theta_j-\theta_i\}.  
$$

In fact, we can choose the following holomorphic  
vector field $\theta_i $ on each open set $\tilde{U}_i$ 
\begin{equation}
\left\{
\begin{array}{l}
	\displaystyle{{\theta }_0=\frac{t-{y_0}^2}{t\ y_0}\frac{\partial }{\partial x_0}+\frac{-2\ x_0\ y_0}{t}\frac{\partial }{\partial y_0}}, \vspace{2mm}\\
	\displaystyle{{\theta }_1=\frac{f_1(x_1,y_1,t)}{t\ F_1(x_1,y_1)}\frac{\partial }{\partial x_1}+\frac{g_1(x_1,y_1,t)}{t\ F_1(x_1,y_1)}\frac{\partial }{\partial y_1}}\vspace{2mm},  \\	
\displaystyle{{\theta }_2=\frac{f_2(x_2,y_2,t)}{t\ F_2(x_2,y_2,t)}\frac{\partial }{\partial x_2}+\frac{g_2(x_2,y_2,t)}{t\ F_2(x_2,y_2,t)}\frac{\partial }{\partial y_2}}.
\end{array}
\right.
\end{equation}
where
\begin{equation}
\renewcommand{\arraystretch}{1.5}
\begin{array}{cl}
f_1(x_1,y_1,t)= & -2 y_1 (t-2 {x_1}^2+5 t x_1 {y_1}^2+9 t {x_1}^2 {y_1}^4+7  t {x_1}^3 {y_1}^6+2 t {x_1}^4 {y_1}^8),  \\
g_1(x_1,y_1,t)= & 1-x_1 {y_1}^2+t {y_1}^4+3 t x_1 {y_1}^6+3 t {x_1}^2 {y_1}^8+t {x_1}^3 {y_1}^{10},  \\
f_2(x_2,y_2,t) = & -t^2+3 t x_2-2 t^3y_2+t x_2 y_2-2 {x_2}^2 y_2+7 t^3 {y_2}^2-8 t^3 {y_2}^3- 8 t^2 x_2 {y_2}^3,  \\
   & +3 t^3 {y_2}^4+18 t^2 x_2 {y_2}^4-10 t^2 x_2 {y_2}^5
   -10 t {x_2}^2 {y_2}^5+11 t {x_2}^2 {y_2}^6-4 {x_2}^3 {y_2}^7,  \\
g_2(x_2,y_2,t)= & -t+t^2 {y_2}^4-2 t^2 {y_2}^5+t^2 {y_2}^6+2 t x_2 {y_2}^6-2 t x_2 {y_2}^7 +{x_2}^2 {y_2}^8.
\end{array}
\end{equation}
Then  we actually have the following  relation as required
$$
\theta_0 = \theta_1, \quad  \theta_2 = \theta_{02} + \theta_0. 
$$

Hence,  we have the differential equation on $ \cS - \cD $ 
as in Theorem \ref{thm:diffeq},  
and on each open set  $ \tilde{U}_i, i=0, 1, 2 $ 
the differential equation can be 
written as follows (cf. (\ref{eq:diffeq})).
\begin{eqnarray}
\mbox{ On $ \tilde{U}_0$ } & & 
\left\{
\begin{array}{l}
	\displaystyle{ \frac{ d x_0}{d t}= - \frac{t - {y_0}^2 }{t y_0} = y_0 \frac{\partial H_0}{\partial y_0} } \\
	\displaystyle{\frac{d y_0}{d t}= \frac{2\ x_0\ y_0}{t}=-y_0\frac{\partial H_0}{\partial x_0} }. 
\end{array} 
\right. \label{eq:diffeqh0} \\ 
\mbox{ On $ \tilde{U}_1 $ } & &  
\left\{\begin{array}{l}
	\displaystyle{ \frac{d x_1}{d t}= - \frac{f_1(x_1,y_1,t)}{t\ F_1(x_1,y_1)} ={F_1(x_1,y_1)}^2 \frac{\partial H_1}{\partial y_1} }  \\
	\displaystyle{ \frac{d y_1}{d t}= - \frac{g_1(x_1,y_1,t)}{t\ F_1(x_1,y_1)}=-{F_1(x_1,y_1)}^2 \frac{\partial H_1}{\partial x_1}. }
\end{array}
\right.    \\ 
\mbox{ On $ \tilde{U}_2 $ } & &
\left\{
\begin{array}{l}
	\displaystyle{ \frac{d x_2}{d t}= - \frac{f_2(x_2,y_2,t)}{t\ F_2(x_2,y_2,t)} } \\
	\displaystyle{ \frac{d y_2}{d t}= - \frac{g_2(x_2,y_2,t)}{t\ F_2(x_2,y_2,t)} } .
\end{array} 
\right. 
\end{eqnarray}

Here $H_0, H_1$ are given by 
\begin{eqnarray}
H_0 & = & \left(- \frac{{x_0}^2}{t}+\frac{y_0}{t}+\frac{1}{y_0}\right) \label{eq:h0},   \\
 H_1 & = & \left({y_1}^2 + x_1 {y_1}^4 + \frac{x_1}{t(1+x_1 {y_1}^2)^2}\right). \label{eq:h1} 
\end{eqnarray}

Moreover on each $\tilde{U}_i$, $i=0, 1, 2$, 
the relative 2-form $\omega_{\cS- \cD}$ are given by 
\begin{eqnarray}
	\omega_{\cS - \cD| \tilde{U}_0 } & = & \frac{1}{y_0}\ dx_0\wedge dy_0 \nonumber\\
	\omega_{\cS - \cD| \tilde{U}_1 }	& = & \frac{1}{{F_1(x_1,y_1)}^2}dx_1\wedge dy_1 \nonumber \\
	\omega_{\cS - \cD| \tilde{U}_2 }& = & \frac{1}{F_2(x_2,y_2,t)}\ dx_2 \wedge dy_2 .\nonumber
\end{eqnarray}
For each $i=0, 1, 2$, consider the 1-form on $\tilde{U_i}$ 
$$
 \theta_i (\omega_{\cS - \cD | \tilde{U}_i }). 
$$
Then since 
$ \theta_i (\omega_{\cS - \cD | \tilde{U}_i })$ 
does not depend on $t$ for $i =0, 1$, 
the fundamental equation (\ref{eq:fund}) is reduced to 
\begin{equation}
d_{\pi} (\theta_i (\omega_{\cS - \cD | \tilde{U}_i })) = 0 , \quad \mbox{ for $ i = 0, 1$ }.
\end{equation}
Though  $\tilde{U}_i$ is not simply connected, we can 
integrate $ \theta_i (\omega_{\cS - \cD | \tilde{U}_i })$ and  obtain $H_i$ for $i = 0, 1$ defined in (\ref{eq:h0}) and (\ref{eq:h1}), that is, 
\begin{equation}
d_{\pi} H_i =  \theta_i (\omega_{\cS - \cD | \tilde{U}_i }).   
\end{equation}

On the other hand, 
since $\theta_2(\omega_{\cS - \cD | \tilde{U}_i})$ is really depend on $t$, the fundamental equation (\ref{eq:fund})  becomes as follows. 
\begin{equation}
\frac{\partial}{\partial t} \left( \frac{1}{F_2(x_2, y_2, t)} \right) dx_2 \wedge d y_2 - d_{\pi} ( \theta_i (\omega_{\cS - \cD | \tilde{U}_i}))= 0. 
\end{equation}
This last equation is equivalent to the following equation, 
which one can check by hand. 
\begin{equation}
\frac{\partial}{\partial t} \left( \frac{1}{F_2} \right) + \frac{\partial}{\partial x_2} \left( \frac{f_2}{t F_2^2} \right) +   \frac{\partial}{\partial y_2} \left( \frac{g_2}{t F_2^2} \right) = 0.  
\end{equation}

Eliminating $x_0$ from the differential equation 
(\ref{eq:diffeqh0}), we obtain the differential equation
\begin{equation}
\frac{d^2 y_0}{d t^2} = \frac{1}{y_0} \left(\frac{d y_0}{d t}\right)^2 - \frac{1}{t}\frac{d y_0}{d t} + \frac{2y_0^2}{t^2} -\frac{2}{t}. 
\end{equation} 
It is easy to see that this equation is equivalent to the equation $P_{III}^{\tilde{D}_8}$ in (\ref{eq:difd8}).

}
\end{Example}

\vspace{1cm}

\begin{center}
{\bf Acknowledgements}
\end{center}

We would like to thank Kazuo Okamoto,  
 Masatoshi Noumi, Hidetaka Sakai, Kyoichi Takano, Akihiro Tuschiya,  Hiroshi Umemura, Yasuhiko Yamada for their continuous  interest in  our work and helpful discussions.  
In particular, H. Umemura sent us  the interesting preprint 
\cite{AL}, H. Sakai sent us his beautiful preprint \cite{Sakai} which have been playing  essential roles in our theory.  Daily  discussions in Kobe  with M. Noumi, K. Takano 
and Y. Yamada have been really stimulating  us to make progress in our research.  
We are really grateful to them. We would like thank Kenji Iohara for helping  us 
to improve our manuscript.  
 
 Thanks are also due to Yuji Shimizu and Kenji Ueno.  They 
organized  the workshops for Painlev\'e equations and Mathematical Physics  both in March 1999 and March 2000, 
which were very stimulating and productive.

\end{document}